\documentclass[11pt]{article}
\usepackage{amscd,amssymb,verbatim,diagrams}
\usepackage[mathscr]{eucal}
\usepackage{mathrsfs}
\usepackage{amsfonts}
\usepackage{amsmath}
\usepackage{amsthm}
\usepackage{latexsym}
\usepackage{eufrak}
\usepackage[all]{xy}

\numberwithin{equation}{section}

\newtheorem{thm}{Theorem}[section]

\newtheorem{lem}[thm]{Lemma}

\newcommand{\Pf}{\noindent {\it Proof}}

\newcommand{\Oo}{\mathcal O}
\newcommand{\Uu}{\mathcal U}

\newcommand{\Pp}{\mathbb P}

\newcommand*{\Dd}{\mathop{\mathrm D\kern0pt}\nolimits}

\long\def\comment#1{}

\renewcommand{\mod}{\operatorname{mod}}

\newcommand{\und}{\underline}

\newcommand{\OO}{{\cal O}}

\newcommand{\DD}{{\cal D}}

\newcommand{\BB}{{\cal B}}

\newcommand{\hra}{\hookrightarrow}

\newcommand{\lan}{\langle}
\newcommand{\ran}{\rangle}

\newcommand{\CC}{{\cal C}}

\newcommand{\Sp}{\operatorname{Sp}}

\renewcommand{\P}{{\Bbb P}}

\newcommand{\ga}{\gamma}
\newcommand{\de}{\delta}
\newcommand{\eps}{\epsilon}

\newcommand{\id}{\operatorname{id}}

\newcommand{\we}{\wedge}

\newcommand{\rk}{\operatorname{rk}}

\renewcommand{\AA}{{\cal A}}

\newcommand{\PP}{{\cal P}}

\newcommand{\Hom}{\operatorname{Hom}}

\newcommand{\Ext}{\operatorname{Ext}}

\renewcommand{\a}{\alpha}
\renewcommand{\b}{\beta}
\newcommand{\om}{\omega}

\newcommand{\la}{\lambda}

\newcommand{\Z}{{\Bbb Z}}

\newcommand{\ot}{\otimes}

\newcommand{\sub}{\subset}
\newcommand{\ed}{\qed\vspace{3mm}}

\newcommand{\GL}{\operatorname{GL}}

\title{Full exceptional collections on the Lagrangian Grassmannians $LG(4,8)$ and $LG(5,10)$}
\author{Alexander Polishchuk and Alexander Samokhin}
\date{}
\begin{document}
\maketitle

\section*{Introduction}

Starting with the seminal works \cite{Beil} and \cite{BGG} on derived categories of coherent sheaves
on projective spaces the techniques involving derived categories have been applied to a variety
of problems in algebraic geometry. Among recent examples one could mention
the relation between semiorthogonal decompositions of derived categories and birational geometry
(see \cite{BO}, \cite{B-flop}), as well as Bridgeland's theory
of stability conditions (see \cite{Bridge}). 
However, there are still some open problems in which not much progress
was made since the 80's. Among them is the problem of describing derived categories of coherent
sheaves on homogeneous varieties. The method of Beilinson in \cite{Beil} was generalized by
Kapranov to the case of quadrics and to partial flag varieties for series $A_n$ (see
\cite{Kap}). Furthermore, it was
realized that the relevant structure is that of a {\it full exceptional collection}, a notion that can be
formulated for an arbitrary triangulated category (see \cite{Rud-sem}). Namely, this is a collection
of objects $E_1, \ldots, E_n$ generating the entire triangulated category with the following vanishing conditions:
$$\Hom^*(E_j,E_i)=0 \text{ for }i<j, \ \Hom^{\neq0}(E_i,E_i)=0,\ \Hom^0(E_i,E_i)=k,$$
where $k$ is the ground field
(which we always assume to be algebraically closed of
characteristic zero). For a smooth projective variety $X$ over $k$ we denote by 
$\Dd^b(X)$ the bounded derived category of coherent sheaves on $X$.
It has been conjectured that for every homogeneous variety $X$ of a semisimple algebraic group
the category $\Dd^b(X)$ admits a full exceptional
collection (of vector bundles).
However, the only homogeneous varieties of simple groups for which this is known 
(other than the examples mentioned above) are:

\noindent  the isotropic Grassmannian of $2$-dimensional planes in a symplectic
$2n$-dimensional space (see \cite{Kuz-iso});

\noindent the isotropic Grassmannian of $2$-dimensional planes in an orthogonal
$2n+1$-dimensional space (see \cite{Kuz-iso});

\noindent the full flag variety for the symplectic and the orthogonal groups (see \cite{Sam});

\noindent the isotropic Grassmannians of a $6$-dimensional symplectic space
(see \cite{Sam});

\noindent the isotropic Grassmannian of $5$-dimensional planes in a $10$-dimensional orthogonal 
space and a certain Grassmannian for type $G_2$ (see \cite{Kuz-hs}).

In the case of the Cayley plane, the minimal homogeneous variety for $E_6$, an exceptional 
collection of $27$ vector bundles, that is conjectured to be full, was constructed in \cite{Manivel}.

In the present paper we construct full exceptional collections of vector bundles in the derived categories
of coherent sheaves of the Lagrangian Grassmannians $LG(4,8)$ and $LG(5,10)$, see
Theorems \ref{exc-col-thm}, \ref{full-thm} and \ref{L5-col-thm}. Note that the situation
is radically different from the previously known cases of classical type in that we have to consider
homogeneous bundles corresponding to reducible representations of the isotropy group.
The new exceptional bundles are constructed as successive extensions of appropriate Schur
functors of the universal quotient bundle.

Checking that the collections we construct are full is done in both cases ``by brute force". 
One needs therefore to find a more conceptual proof before
trying to generalize our results to other Lagrangian Grassmannians.
It seems plausible that an exceptional collection $(E_1,\ldots,E_n)$ in $\Dd^b(X)$ such that
classes of $E_i$ generate the Grothendieck group $K_0(X)$, is automatically full.
Thus, it would be enough to have $n=\rk K_0(X)$.
Recall that a full triangulated subcategory $\CC\sub\DD$ generated by an exceptional collection
in $\DD$ is {\it admissible} (see \cite{Bondal}, Thm. 3.2). By definition, this means that the
inclusion functor $\CC\to\DD$ admits left and right adjoint functors $\DD\to\CC$. 
To check that $\CC=\DD$ is equivalent to showing that the right orthogonal $\CC^{\perp}\sub\DD$
is zero, where $\CC^{\perp}=\{A\in\DD\ |\ \Hom_{\DD}(\CC,A)=0\}$. It is known that
$\CC^{\perp}$ is also admissible. Thus, the above statement 
would follow from the Nonvanishing conjecture of A.~Kuznetsov (see \cite{Kuz-Hoch}, Conjecture 9.1 and Corollary 9.3) that a nonzero admissible subcategory should have nonzero Hochschild homology.

\vspace{2mm}

\noindent
{\it Acknowledgments}.
The work of the first author was partially supported by the NSF grant DMS-0601034.
Part of this work was done while the second author was visiting the University of Oregon and
the IHES.
He gratefully acknowledges the hospitality and support of both institutions.

\section{Applications of the Bott's theorem in the case of Lagrangian Grassmannians}

Let $V$ be a symplectic vector space of dimension $2n$.
Consider the Largangian Grassmannian $LG(V)$ of  $V$
(we also use the notation $LG(n,2n)$).
We have the basic exact sequence of vector bundles on $LG(V)$
\begin{equation}\label{basic-seq}
0\to \Uu\to V\ot\OO \to Q\to 0
\end{equation}
where $\Uu=Q^{\ast}$ is the tautological subbundle, and $Q$ is the tautological quotient-bundle.
We set $\OO(1)=\wedge^n Q$. This is an ample generator of the Picard group of $LG(V)$.
It is well known that the canonical line bundle on $LG(V)$ is isomorphic to $\OO(-n-1)$.

The variety $LG(V)$ is a homogeneous space for the symplectic group $\Sp(V)=\Sp(2n)$. Namely,
it can be identified with $\Sp(2n)/P$, where $P$ is the maximal parabolic associated with the
simple root $\a_n$. Here we use the standard numbering of the vertices in the Dynkin diagram 
$D_n$ as in \cite{Bour}. Recall that the semisimple part of $P$ is naturally identified with
$\GL(n)$. Thus, to every representation of $\GL(n)$ one can associate a homogeneous vector
bundle on $LG(V)$. This correspondence is compatible with tensor products and the standard
representation of $\GL(n)$ corresponds to $Q$. For our purposes it will be convenient to
identify the maximal torus of $\Sp(2n)$ with that of $\GL(n)\sub P$. One can easily check
that under this identification the half-sum of all the positive roots of $\Sp(2n)$ is equal to
$$\rho=n\eps_1+(n-1)\eps_2+\ldots+\eps_n,$$
where $(\eps_i)$ is the standard basis of the weight lattice corresponding to $\GL(n)$.
Note that with respect to this basis the roots of $\Sp(2n)$ are $\pm\eps_i$ and
$\pm\eps_i\pm\eps_j$. Thus, a weight $x_1\eps_1+\ldots+x_n\eps_n$ is singular for $\Sp(2n)$
if and only if
either there exists $i$ such that $x_i=0$, or there exist $i\neq j$ such that $x_i=\pm x_j$.
The Weyl group $W$ of $\Sp(2n)$ is the semidirect product of $S_n$ and $\Z_2^n$
acting by permutations and sign changes $x_i\mapsto -x_i$.
A weight $x_1\eps_1+\ldots+x_n\eps_n$ is dominant for $\Sp(2n)$ if and only if 
$x_1\ge x_2\ge \ldots\ge x_n\ge 0$.

For a dominant weight $\lambda=(a_1,\ldots,a_n)$ of $\GL(n)$ (where $a_1\ge a_2\ge\ldots\ge a_n$), 
let $S^{\lambda}$ denote the corresponding Schur functor (sometimes we omit the tail of
zeros in $\lambda$). Note that by definition, $S^{(a_1+1,\ldots,a_n+1)}=\det\ot S^{(a_1,\ldots,a_n)}$. Hence,
$$S^{(a_1+1,\ldots,a_n+1)}Q\simeq S^{(a_1,\ldots,a_n)}Q(1).$$

Our main computational tool is Bott's theorem on cohomology of homogeneous vector bundles.
In the case of the Lagrangian Grassmannian $LG(V)$ it states the following.

\begin{thm} (Theorem $IV'$ of \cite{Bott}) 
\begin{enumerate}
\item If $\la+\rho$ is singular then $H^{\ast}(LG(V),S^{\lambda}Q)=0$;
\item if $\la+\rho$ is non-singular and $w\in W$ is an element of minimal length $\ell$ such that
$\mu=w(\la+\rho)-\rho$ is dominant for $\Sp(2n)$, then $H^i(LG(V),S^{\lambda}Q)=0$ for $i\neq\ell$ and
$H^{\ell}(LG(V),S^{\la}Q)$ is an irreducible representation of $\Sp(2n)$ with the highest weight $\mu$.
\end{enumerate}
\end{thm}

Below we will often abbreviate 
$H^{\ast}(LG(V),?)$ to $H^{\ast}(?)$.

\begin{lem}\label{coh-lem} 
One has

\noindent
(i) $H^{\ast}(\OO(i))=0$ for $i\in [-n,-1]$; $H^{>0}(\OO)=0$ and $H^0(\OO)=k$.

\noindent
(ii) $H^{\ast}(\wedge^k Q(i))=0$ for $k\in [1,n-1]$ and $i\in [-n-1,-1]$. Also,
for $k\in [1,n-1]$ one has $H^{>0}(\wedge^k Q)=0$ and
$H^0(\wedge^k Q)$ is an irreducible representation of $\Sp(2n)$ with the highest weight
$((1)^k,(0)^{n-k})$ ($k$ $1$'s).
\end{lem}

\Pf . (i) We have in this case $\la+\rho=(n+i,\ldots,1+i)$ which is singular fo $i\in [-n,-1]$.
For $i=0$ we have $\la+\rho=\rho$.

\noindent
(ii) The bundle $\wedge^k Q$ corresponds to the weight $((1)^k,(0)^{n-k})$.
Thus, $\wedge^k Q(i)$ corresponds to $\la=((1+i)^k,(i)^{n-k})$, so
$\la+\rho=(n+1+i,\ldots, n-k+2+i, n-k+i,\ldots,1+i)$. In the case when $i\in [-n-1, -n-2+k]$ or
$i\in [-n+k,-1]$ one of the coordinates is zero. On the other hand, for $i=-n-1+k$ the
sum of the $k$th and $(k+1)$st coordinates is zero. Hence, $\la+\rho$ is singular for $i\in [-n-1,-1]$.
\ed

When computing the $\Ext$-groups on $LG(V)$ between the bundles of the form $S^{\lambda}Q$
it is useful to observe that
$$(S^{(a_1,\ldots,a_n)}Q)^{\ast}\simeq S^{(a_1-a_n,a_1-a_{n-1},\ldots,0)}(-a_1).$$
To compute the tensor products of the Schur functors we use Littlewood-Richardson rule.

\begin{lem}\label{gen-exc-lem}  Assume that $n\ge 3$.

\noindent (i)
One has $\Hom^*(\wedge^k Q,\wedge^l Q(i))=0$ for $i\in [-n,-1]$ and $k,l\in [0,n-2]$.
Also, $\Hom^*(\wedge^k Q,\wedge^l Q)=0$ for $k,l\in [0,n-2]$ and $k>l$. 
All the bundles $\wedge^k Q$ are exceptional.

\noindent (ii)
For $k<n$ one has $\Hom^*(\wedge^k Q,\wedge^{k+1}Q)=V$ (concentrated in degree $0$).
Furthermore, the natural map $Q\to\und{\Hom}(\wedge^k Q,\wedge^{k+1}Q)$ induces an
isomorphism on $H^0$.
\end{lem}

\Pf . (i) Recall that
$$\wedge^k Q^{\ast}=\wedge^{n-k} Q(-1)=S^{((1)^{n-k},(0)^k)}Q(-1).$$
Therefore, for $k>l$, $k+l\neq n$, the tensor product
$\wedge^k Q^{\ast}\ot \wedge^l Q\simeq \wedge^{n-k}Q\ot \wedge^l Q(-1)$ decomposes into
direct summands of the form $S^{\la}Q$ with $\la=((1)^a,(0)^b,(-1)^c)$, where $b>0$ and $c>0$.
It is easy to see that in this case $\la+((i)^n)+\rho$ will be singular for $i\in [-n,0]$. Furthermore, 
even if $k+l=n$ but $l<k<n-1$, we claim that the weights $\la+((i)^n)+\rho$ will still be singular
for $i\in[-n,0]$. Indeed, this follows easily from the fact that $\la=((1)^a,(0)^b,(-1)^c)$ with $c>0$, and
either $b>0$ or $c>1$ or $a>1$. 
Hence, $\Hom^*(\wedge^k Q,\wedge^lQ(i))=0$, where $i\in [-n,0]$, $n>k>l\ge 0$ and $(k,l)\neq (n-1,1)$.
Using Serre duality we deduce the needed vanishing for the case $k<l$.
In the case when $k=l$ the tensor product
$\wedge^k Q^{\ast}\ot\wedge^k Q\simeq \wedge^{n-k}Q\ot\wedge^k Q(-1)$ will contain exactly on summand isomorphic
to $\OO$, and the other summands of the same form as above with $c>0$. The same argument as
before shows that $\Hom^*(\wedge^kQ,\wedge^kQ(i))=0$ for $i\in [-n,-1]$ and that
$\Hom^*(\wedge^kQ,\wedge^kQ)=k$ (concentrated in degree $0$).

\noindent
(ii) The tensor product $\wedge^k Q^{\ast}\ot\wedge^{k+1}Q\simeq\wedge^{n-k}Q\ot\wedge^{k+1}Q(-1)$
decomposes into the direct sum of $Q$ and of summands
of the form $S^{\la}Q$ with $\la=((1)^a,(0)^b,(-1)^c)$, where $c>0$. In the latter case
the weight $\la+\rho$ is singular, so these summands do not contribute to cohomology.
\ed

Next, for $k\in [1,n-3]$
consider the vector bundle $R_k:=S^{(2,(1)^k)}Q$, so that we have a direct sum decomposition
$$Q\ot \wedge^{k+1} Q= \wedge^{k+2} Q\oplus R_k.$$
One can check that $R_k$ itself is not exceptional but in the next section
we are going to construct a related exceptional bundle on $LG(V)$.

\begin{lem}\label{R-van-lem}
For $1\le k\le n-3$, $0\le l\le n-2$ and $-n\le i\le -1$ one has
$$\Hom^*(\we^l Q, R_k(i))=\Hom^*(R_k,\we^l Q(i))=0.$$ 
Furthermore, for $l>k+1$ one has $\Hom^*(\we^l Q, R_k)=0$, while
for $l<k$ one has $\Hom^*(R_k,\we^l Q)=0$.
\end{lem}

\Pf .
By Littlewood-Richardson rule, the tensor product 
$\we^l Q^*\ot R_k=\we^{n-l} Q\ot R_k(-1)$ decomposes into direct summands
of the form $S^{\la}$, where $\la$ has one of the following types:

\noindent
(i) $\la=(1,(0)^{k+n-l},(-1)^{l-k-1})$, provided $l\ge k+1$ (note that $k+n-l\ge 3$);

\noindent
(ii) $\la=((1)^a,(0)^b,(-1)^c)$, where $1\le a\le k+1$, $a+b\ge k+1$, $a+b+c=n$, $2a+b=k+n-l+2$;

\noindent
(iii) $\la=(2,(1)^a,(0)^b,(-1)^c)$, where $a\le k$, $a+b\ge k$, $a+b+c=n-1$, $2a+b=k+n-l-1$.

\noindent
In case (i) the weight $\la+((i)^n)+\rho$ will be singular for $i\in [-n-1,-1]$. In the case
$l>k+1$ it will also be singular for $i=0$. Next, let us consider case (ii).
If $b>0$ then the weight $\la+((i)^n)+\rho$ will be singular for $i\in [-n-1,-1]$ and if in addition
$c>0$ then it will be also singular for $i=0$. Note that the case $c=0$ occurs only when $l\le k+1$.
In the case $b=0$ we should have $a=k+1$, so $2\le a\le n-2$, which implies that $\la+((i)^n)+\rho$ is 
singular for $i\in [-n-1,0]$. Finally, let us consider case (iii).
If $a>0$, $b>0$ and $c>0$ then the weight $\la+((i)^n)+\rho$ will be singular for $i\in [-n-1,0]$.
The case $c=0$ can occur only when $l\le k$. In the case $b=0$ we should have $a=k$,
so $c=n-k-1\ge 2$ which implies that the above weight is still singular for $i\in [-n-1,0]$.
In the case $a=0$ we have $b=k+n-l-1\ge 2$, so we deduce that the above weight will be singular
for $i\in [-n,0]$. Note that the case $a=0$ can occur only for $l\ge k$.

The above analysis shows the vanishing of $\Hom^*(\we^l Q, R_k(i))$ for $i\in [-n,-1]$, as well
as vanishing of $\Hom^*(\we^l Q, R_k)$ for $l>k+1$ and of $\Hom^*(\we^l Q, R_k(-n-1))$ for
$l<k$. Applying Serre duality we deduce the remaining assertions. 
\ed

Note that in the above lemma we have skipped the calculation of $\Hom^*(R_k,\we^l Q)$ and
$\Hom(\we^l Q,R_k)$ for $l=k$ and $l=k+1$. This will be done in the following lemma,
where we also prove a number of other auxiliary statements.
Let us consider a natural map $f:V\ot\wedge^{k+1}Q\to R_k$ induced by the projection
$Q\ot\wedge^{k+1} Q\to R$ and the map $V\ot \OO\to Q$.

\begin{lem}\label{Bott-lem} 
Assume $1\le k\le n-3$. Then one has 

\noindent
(i) $\Hom^{>0}(\wedge^{k+1} Q,R_k)=0$ and $\Hom^0(\wedge^{k+1} Q,R_k)=V$.
The map $f$ induces an isomorphism on $\Hom^*(\wedge^{k+1} Q, ?)$.

\noindent
(ii) One has $\Hom^{>0}(\wedge^k Q,R_k)=0$. Also,
the natural map $Q\ot Q\to\und{\Hom}(\wedge^k Q,R_k)$ induces an isomorphism on $H^0$, so that
$\Hom^0(Q,R_k)\simeq V\ot V/k$.

\noindent
(iii) $\Hom^1(R_k,\wedge^k Q)=k$, 
$\Hom^1(R_k,\wedge^{k+1} Q)=V$, $\Hom^{\neq 1}(R_k,\wedge^k Q)=
\Hom^{\neq 1}(R_k,\wedge^{k+1} Q)=0$.
The natural map $S^2Q^{\ast}\to\und{\Hom}(R_k,\wedge^k Q)$ induces an isomorphism on $H^1$.

\noindent
(iv) $\Hom^{>1}(R_k,R_k)=0$, $\Hom^0(R_k,R_k)=k$, $\Hom^1(R_k,R_k)=V\ot V/k$.

\noindent
(v) $H^{\ast}(Q^{\ast}\ot S^2Q^{\ast})=0$.

\noindent
(vi) $H^i(Q^{\ast}\ot Q\ot S^2Q^{\ast})=0$ for $i\neq 1$.
\end{lem}

\Pf . (i) By Littlewood-Richardson rule we have 
$$\wedge^{k+1} Q^{\ast}\ot R_k\simeq 
\wedge^{n-k-1} Q(-1)\ot S^{(2,(1)^k)}Q\simeq
Q\oplus \ldots
$$
where the remaining summands correspond to highest weights
$\la=(a_1,\ldots,a_n)$ such that $a_n=-1$.
For such $\la$ the weight $\la+\rho$ is singular,
hence these summands do not contribute to cohomology. Thus, the unique embedding of
$Q$ into $\und{\Hom}(\wedge^{k+1} Q,R_k)$ induces an isomorphism on cohomology.
This immediately implies the result (recall that $H^{\ast}(Q)=V$ by Lemma \ref{coh-lem}).

\noindent (ii) Applying Littlewood-Richardson rule again we find
$$\wedge^k Q^{\ast}\ot R_k\simeq S^2Q\oplus\wedge^2Q\oplus\ldots$$
where the remaining summands correspond to highest weights $\la=(a_1,\ldots,a_n)$
with $a_n=-1$. The sum of the first two terms is exactly
the image of the natural embedding $Q\ot Q\to\und{\Hom}(Q,R)$.

\noindent
(iii) We have 
$$R_k^{\ast}\ot \wedge^k Q\simeq S^2Q^{\ast}\oplus\ldots$$
where all the remaining summands correspond to highest weights $\la=(a_1,\ldots,a_n)$
such that either $a_n=-1$ or $(a_{n-1},a_n)=(-1,-2)$. In both cases $\la+\rho$ is singular,
hence these summands do not contribute to cohomology.
for $S^2Q^{\ast}=S^{((0)^{n-1},-2)}Q$ one has
$\la+\rho=(n,\ldots,2,-1)$. Hence, applying a simple reflection
we get exactly $\rho$. This means that only $H^1$ is nonzero, and it is one-dimensional.

Similarly,
$$R_k^{\ast}\ot \wedge^{k+1}Q\simeq S^{(1,(0)^{n-2},-2)}Q\oplus\ldots$$
where the remaining summands have singular $\la+\rho$.
For $\la=(1,(0)^{n-2},-2)$ we have
$\la+\rho=(n+1,\ldots,2,-1)$. This differs by a single reflection from $\rho+(1,(0)^{n-1})$.
Hence only $H^1$ is nonzero and $H^1(R_k^{\ast}\ot\wedge^{k+1}Q)\simeq V$.

\noindent (iv) We have
$$R_k^{\ast}\ot R_k\simeq S^{((2)^{n-2},1,0)}Q(-2)\ot S^{(2,1)}Q\simeq
S^{(2,(0)^{n-2},-2)}Q\oplus S^{(1,1,(0)^{n-3},-2)}Q\oplus \OO\oplus\ldots,$$
where the remaining terms do not contribute to cohomology.
The first two terms contribute only to $H^1$. Namely, the corresponding weights
$\la+\rho$ differ by a single reflection from $\rho+(2,(0)^{n-1})$ and $\rho+(1,1,(0)^{n-2})$,
respectively.

\noindent (v)
We have 
$$Q^{\ast}\ot S^2Q^{\ast}\simeq S^3Q^{\ast}\oplus S^{(2,1)}Q^{\ast}\simeq 
S^{((0)^{n-1},3)}Q\oplus S^{((0)^{n-2},-1,-2)}Q.$$
In both cases $\la+\rho$ is singular.

\noindent (vi)
We have 
$$Q\ot Q^{\ast}\ot S^2Q^{\ast}=Q\ot S^3Q^{\ast}\oplus Q\ot S^{(2,1)}Q^{\ast}
\simeq (S^{((0)^{n-1},-2)}Q)^{\oplus 2}\oplus\ldots,
$$
where the remaining summands do not contribute to cohomology.
For the first summand we have $\la+\rho=(n,\ldots,2,-1)$ which is obtained by applying 
a simple reflection to a dominant weight. Hence, the cohomology is concentrated in degree $1$.
\ed

\section{A family of exceptional vector bundles on $LG(V)$}

Let us fix $k\in [1,n-3]$.
The natural map $f:V\ot\wedge^{k+1} Q\to Q\ot\wedge^{k+1} Q$ is surjective, so we 
obtain an exact sequence of vector bundles
\begin{equation}\label{SR-seq}
0\to S_k\to V\ot\wedge^{k+1} Q\stackrel{f}{\to} R_k\to 0.
\end{equation}
Using the composite nature of $f$ we also get an exact sequence
\begin{equation}\label{S-seq}
0\to Q^{\ast}\ot\wedge^{k+1} Q\to S_k\to \wedge^{k+2} Q\to 0.
\end{equation}
We have a natural embedding of vector bundles
$$\wedge^k Q\hra \und{\Hom}(Q,\wedge^{k+1} Q)=Q^{\ast}\ot\wedge^{k+1} Q\hra S_k.$$
Now we define $E_k$ to be the quotient $S_k/\wedge^k Q$, so that we have an exact sequence
\begin{equation}\label{QSE-seq}
0\to \wedge^k Q\to S_k\to E_k\to 0.
\end{equation}

\begin{lem}\label{split-lem} 
The exact sequence \eqref{QSE-seq} splits canonically, so we have
$S_k\simeq\wedge^k Q\oplus E_k$. Furthermore, the bundles
$\wedge^k Q$ and $E_k$ are orthogonal to each other, i.e.,
$$\Hom^*(\wedge^k Q,E_k)=\Hom^*(E_k,\wedge^k Q)=0.$$
\end{lem}

\Pf . First, we claim that $\Hom^0(S_k,\wedge^k Q)=k$ and $\Hom^i(S_k,\wedge^k Q)=0$
for $i\neq 0$.
Indeed, this follows immediately from the exact sequence \eqref{SR-seq} and from 
Lemma \ref{Bott-lem}(iii) since $\Hom^*(\wedge^{k+1} Q,\wedge^k Q)=0$
by Lemma \ref{gen-exc-lem}. Next, using the vanishing of $\Hom^*(\wedge^{k+2} Q,\wedge^k Q)$
and the exact sequence \eqref{S-seq} we see that the embedding
$Q^{\ast}\ot\wedge^{k+1} Q\hra S_k$ induces an isomorphism
on $\Hom^*(?,\wedge^k Q)$. Hence, the nonzero morphism $S_k\to\wedge^k Q$ restricts
to the nonzero morphism $Q^{\ast}\ot\wedge^{k+1} Q\to\wedge^k Q$, unique up to scalar. 
The latter morphism is proportional to the natural contraction operation. Hence,
its restriction to $\wedge^k Q\sub Q^{\ast}\ot\wedge^{k+1} Q$ is nonzero. Therefore, we get
a splitting of \eqref{QSE-seq}.
The vanishing of $\Hom^*(E_k,\wedge^k Q)$ also follows. On the other hand, from the exact
sequence \eqref{SR-seq}, using Lemma 
\ref{Bott-lem}(ii) we get $\Hom^0(\wedge^k Q,S_k)=k$ and $\Hom^i(\wedge^k Q,S_k)=0$ for
$i\neq 0$. This implies that $\Hom^*(\wedge^k Q, E_k)=0$. 
\ed

By the above lemma we have a unique morphism $S_k\to\wedge^k Q$ extends the identity
morphism from $\wedge^k Q\sub S_k$. Pushing forward the extension given by \eqref{SR-seq}
under this morphism we get an extension
\begin{equation}\label{FR-seq}
0\to \we^k Q\to F_k\to R_k\to 0.
\end{equation}
Furthermore, we also get an exact sequence
\begin{equation}\label{EF-seq}
0\to E_k\to V\ot \we^{k+1}Q\to F_k\to 0.
\end{equation}

Let us recall the definition of the mutation operation.
For an exceptional pair $(A,B)$ in a triangulated category $\DD$,
the right mutation is a pair $(B,R_{B}A)$, where
  $R_{B}A$ is defined to be a cone of the triangle
\begin{equation}\nonumber
\dots \longrightarrow R_{B}A[-1] \longrightarrow A\longrightarrow \Hom _{\DD}^{\bullet}(A,B)^{\ast}\otimes B\longrightarrow
  R_{B}A\longrightarrow \dots \quad .
\end{equation}
The pair $(B,R_B A)$ is again exceptional. 

\begin{thm}\label{exc-thm} Let $k\in [1,n-3]$. 
The bundle $F_k$ is the unique nontrivial extension of $R_k$ by $\we^k Q$.
The bundles $E_k$ and $F_k$ are exceptional, and $F_k$ is the right mutation of
$E_k$ through $\we^{k+1}Q$. Also, one has $F_k^{\ast}(1)\simeq E_{n-2-k}$.
\end{thm}

\Pf . {\bf Step 1}. $\Hom^*(\we^{k+1}Q,E_k)=\Hom^*(F_k,\we^k Q)=0$. Indeed, the first vanishing
follows immediately from the exact sequence \eqref{SR-seq}. The second vanishing
follows from the exact sequence  
\eqref{EF-seq} since $\Hom^*(\we^{k+1}Q,\we^k Q)=0$ and $\Hom^*(E_k,\we^k Q)=0$ by Lemma
\ref{split-lem}.

 \noindent
{\bf Step 2}. $F_k$ is a nontrivial extension of $R_k$ by $\we^k Q$ (recall that by Lemma
\ref{Bott-lem}(iii) there is a unique such extension). Indeed, 
otherwise we would have a surjective map $F_k\to\we^k Q$ which is impossible by Step 1.

\noindent
{\bf Step 3}. $E_k$ is isomorphic to $F_{n-2-k}^{\ast}(1)$.
We have $Q^{\ast}\ot\we^{k+1}Q\simeq \we^k Q\oplus R^{\ast}_{n-k-2}(1)$. Therefore, from
the exact sequence \eqref{S-seq} we get an exact sequence
$$0\to R^{\ast}_{n-2-k}(1)\to E_k\to \we^{k+2}Q\to 0.$$
We claim that it does not split. Indeed, otherwise we would get an inclusion
$\we^{k+2}Q\hra E_k$ which is impossible since $\Hom(\we^{k+1}Q,\we^{k+2}Q)\neq 0$
but $\Hom(\we^{k+1}Q,E_k)=0$. Comparing this with the extension \eqref{FR-seq} for $n-2-k$
instead of $k$ we get the result.

\noindent
{\bf Step 4}. The natural map
$$H^0(Q\ot Q)\ot H^1(S^2 Q^{\ast})\to H^1(Q\ot Q\ot S^2 Q^{\ast})$$
is an isomorphism. Indeed, it is easy to check using Bott's theorem that both sides are isomorphic to 
$V^{\ot 2}/k$, so it is enough to check surjectivity. Therefore, it suffices to check 
surjectivity of the maps 
$$H^0(Q)\ot H^1(S^2 Q^{\ast})\to H^1(Q\ot S^2 Q^{\ast}) \ \text{and}$$
$$H^0(Q)\ot H^1(Q\ot S^2 Q^{\ast})\to H^1(Q\ot Q\ot S^2 Q^{\ast}).$$
Using the exact sequence \eqref{basic-seq} we deduce this from the vanishing
of $H^2(Q^{\ast}\ot S^2 Q^{\ast})$ and $H^2(Q^{\ast}\ot Q\ot S^2 Q^{\ast})$
(see Lemma \ref{Bott-lem}(v),(vi)).

\noindent
{\bf Step 5}. The composition map
$$\Hom^0(\we^k Q,R_k)\otimes\Hom^1(R_k,\we^k Q)\to\Hom^1(R_k,R_k)$$
is an isomorphism. Note that by Lemma \ref{Bott-lem}(ii),(iii),(iv),
both sides are isomorphic to $V^{\ot 2}/k=S^2V\oplus \we^2 V/k$, so it is enough
to check surjectivity. Let us define the natural morphisms
$$\a:S^2Q^{\ast}\to R_k^{\ast}\ot \we^k Q,$$
$$\b:Q\ot Q\to \we^k Q^{\ast}\ot R_k,$$
as follows. Consider the Koszul complex for the symmetric algebra $S^*Q$
$$0\to\we^{k+2}Q\stackrel{d_1}{\to}Q\ot\we^{k+1} Q\stackrel{d_2}{\to} S^2Q\ot\we^k Q\to\ldots$$
Then $R_k$ can be identified with the image of $d_2$ (or cokernel of $d_1$). 
In particular, we have a natural
embedding $R_k\to S^2Q\ot\we^k Q$ which induces $\a$ by duality.
On the other hand, the natural projection $Q\ot\we^{k+1} Q\to R_k$ gives rise to
the composed map 
$$Q\ot Q\ot\we^k Q\stackrel{\id_Q\ot\mu_k}{\to} Q\ot\we^{k+1}Q\to R_k$$
where $\mu_k:Q\ot\we^k Q\to\we^{k+1}Q$ is given by the exterior product. The map $\b$
is obtained from the above map by duality.
The morphisms $\a$ and $\b$ can be combined into a map
$$\ga:S^2Q^{\ast}\ot Q\ot Q\stackrel{\a\ot\b}{\to}R_k^{\ast}\ot \we^k Q\ot \we^k Q^{\ast}\ot R_k\to
R_k^{\ast}\ot R_k,$$
where the last arrow is induced by the trace map on $\we^k Q$.
By Step 4, it remains to check that the maps $\a$, $\b$ and $\ga$ induce isomorphisms on
cohomology. In fact, we are going to prove that all these maps are embeddings of a direct summand
by constructing the maps $p_{\a}$, $p_{\b}$ and $p_{\ga}$ in the opposite direction such that
$p_{\a}\circ\a$, $p_{\b}\circ\b$ and $p_{\ga}\circ\ga$ are proportional to identity. To this end we use the Koszul complex for the exterior algebra $\we^*Q$
$$\ldots\to S^2 Q\ot\we^k Q\stackrel{\de_2}{\to}Q\ot\we^{k+1}Q\stackrel{\de_1}{\to}\we^{k+2}Q\to 0$$
We can identify $R_k$ with the kernel of $\de_1$ (or image of $\de_2$). Hence,
we have natural map $S^2 Q\ot\we^k Q\to R_k$. By duality this corresponds to a map
$$p_{\a}:R_k^{\ast}\ot\we^k Q\to S^2Q^{\ast}.$$
On the other hand, we have a natural embedding
$$R_k\to Q\ot\we^{k+1}Q\to Q\ot Q\ot\we^k Q$$
that gives rise to a map $p_{\b}:\we^k Q^{\ast}\ot R_k\to Q\ot Q$.
Combining $p_{\a}$ and $p_{\b}$ we obtain a map
$$p_{\ga}:R_k^{\ast}\ot R_k\to R_k^{\ast}\ot\we^k Q\ot\we^kQ^{\ast}\ot R_k\to S^2Q^{\ast}\ot Q\ot Q.$$
A routine calculation proves our claim about the compositions $p_{\a}\circ\a$, $p_{\b}\circ\b$ and
$p_{\ga}\circ\ga$.

\noindent
{\bf Step 6}. Now we can prove that $F_k$ is exceptional (and hence, $E_k$ is also exceptional
by Step 3).
Applying the functor $\Hom^*(F_k,?)$ to the exact sequence 
\eqref{FR-seq} and using Step 1 we get isomorphisms $\Hom^i(F_k,F_k)\simeq \Hom^i(F_k,R_k)$.
Next, applying the functor $\Hom^*(?.R_k)$ to the same sequence we get a long exact sequence
$$\ldots\to\Hom^{i-1}(R_k,\we^k Q)\to\Hom^i(R_k,R_k)\to \Hom^i(F_k,R_k)\to\Hom^i(R_k,\we^k Q)\to\ldots$$
It remains to apply Lemma \ref{Bott-lem}(iii) and Step 5 to conclude that $\Hom^i(F_k,R_k)=0$ for $i>0$ and
$\Hom^0(F_k,R_k)=k$.

\noindent
{\bf Step 7}. To check that $F_k$ is the right mutation of $E_k$ through $\we^{k+1}Q$ it remains
to prove that $\Hom^i(E_k,\we^{k+1}Q)=0$ for $i\neq 0$ and $\Hom^0(E_k,\we^{k+1}Q)$.
Applying the functor $\Hom^*(?,\we^{k+1}Q)$ to the sequence \eqref{SR-seq} we get
by Lemma \ref{Bott-lem}(iii) an exact sequence
$$0\to V\to \Hom^0(S_k,\we^{k+1}Q)\to V\to 0$$
along with the vanishing of $\Hom^{>0}(S_k,\we^{k+1}Q)$. Since $S_k=\we^k Q\oplus E_k$,
the assertion follows.
\ed

We are going to compute some $\Hom$-spaces involving the bundles $E_k$ that we will need
later.

\begin{lem}\label{E-van-lem} 
Assume that $l\in [0,n-2]$ and $k\in [1,n-3]$. Then
for $i\in [-n,-1]$ one has
$$\Hom^*(\we^l Q, E_k(i))=\Hom^*(E_k,\we^l Q(i))=0.$$
For $l>k$ one has $\Hom^*(\we^l Q,E_k)=0$, while for $l<k$ one has
$\Hom(E_k,\we^l Q)=0$ (recall that for $l=k$ both these spaces vanish by
Lemma \ref{split-lem}).
\end{lem}

\Pf  . It is enough to check similar assertions with $S_k$ instead of $E_k$.
Using the exact sequence \eqref{SR-seq} we reduce the required vanishing for $i\in [-1,-n]$
to Lemmas \ref{gen-exc-lem}(i) and \ref{R-van-lem}. To prove the remaining vanishings
we use in addition the fact that $\Hom^*(\we^{k+1}Q, S_k)=0$ that 
follows from Lemma \ref{Bott-lem}(i).
\ed

\section{The case of $LG(4,8)$}

Now let us assume that $V$ is $8$-dimensional. Let $E=E_1$.




\begin{thm}\label{exc-col-thm} 
The following collection on $LG(4,8)$ is exceptional:
$$(\OO,E,Q,\wedge^2Q,\OO(1),Q(1),\wedge^2(1),\ldots,\OO(4),Q(4),\wedge^2Q(4)).$$
\end{thm}

\Pf . We already know that all these bundles are exceptional.
The required orthogonality conditions follow from
Lemma \ref{gen-exc-lem}, Lemma \ref{split-lem}, and Lemma \ref{E-van-lem}.
\ed

\begin{lem}\label{full-lem}
Let $\CC\sub \Dd^b(LG(4,8))$ be the triangulated subcategory generated by the exceptional
collection in Theorem \ref{exc-col-thm}. Then the following bundles belong to $\CC$:

\noindent
(i) $Q^{\ast}(j)$, $j=0,\ldots,4$;

\noindent
(ii) $S^2Q(j)$, $j=0,\ldots,4$;

\noindent
(iii) $Q\ot\wedge^2 Q(j)$, $j=0,\ldots,3$;

\noindent
(iv) $Q\ot Q^{\ast}(j)$, $j=1,\ldots,4$.
\end{lem}

\Pf . {\bf Step 1}. $Q^{\ast}(j),S^2Q^{\ast}(j)\in\CC$ for $j=0,\ldots,4$. 
Indeed, the fact that $Q^{\ast}(j)\in\CC$ follows immediately from \eqref{basic-seq}. Similarly,  
the assertion for $S^2Q^{\ast}(j)$ follows from the exact sequence
\begin{equation}\label{S2Q*-seq}
0\to S^2Q^{\ast}\to S^2V\ot\OO\to V\ot Q\to\wedge^2Q\to 0
\end{equation}
obtained from \eqref{basic-seq}.

\vspace*{0.1cm}

\noindent
{\bf Step 2}. $Q\ot Q(j)\in\CC$ for $j=1,2,3,4$. This follows from the exact sequence
\begin{equation}\label{S2Q-seq}
0\to \wedge^2 Q^{\ast}\to V\ot Q^{\ast}\to S^2V\ot\OO\to S^2Q\to 0,
\end{equation}
dual to \eqref{S2Q*-seq}, since $\we^2Q^{\ast}=\we^2Q(-1)$ and $Q^{\ast}(j)\in\CC$ by Step 1.

\vspace*{0.1cm}

\noindent
{\bf Step 3}. $\wedge^2Q\ot\wedge^2Q(2)\in\CC$.
It follows from the basic sequence \eqref{basic-seq} that $\wedge^4 V\ot\OO(3)$ has
a filtration with the subsequent quotients $\OO(4)$, $Q^{\ast}\ot Q^{\ast}(4)$,
$\wedge^2Q\ot\wedge^2Q(2)$, $Q\ot Q(2)$ and $\OO(2)$. All of them except for 
$\wedge^2Q\ot\wedge^2Q(2)$ belong to $\CC$, by Steps 1 and 2. This
implies the assertion.

\vspace*{0.1cm}

\noindent
{\bf Step 4}. $Q^{\ast}\ot Q(j)\in\CC$ for $j=1,2,3,4$. Indeed, tensoring \eqref{basic-seq} with $Q$
we get an exact sequence
$$0\to Q^{\ast}\ot Q\to V\ot Q\to Q\ot Q\to 0,$$
so the assertion follows from Step 2.

\vspace*{0.1cm}

\noindent
{\bf Step 5}. $Q\ot\wedge^2Q\in\CC$. First, observe that $S_1=Q\oplus E\in\CC$.
Now the exact sequence \eqref{SR-seq} shows that $R_1\in\CC$.
But $Q\ot\wedge^2Q=\wedge^3Q\oplus S^{(2,1)}Q=Q^{\ast}(1)\oplus R_1$, so it is in $\CC$ (recall
that $Q^{\ast}(1)\in\CC$ by Step 1).

\vspace*{0.1cm}

\noindent
{\bf Step 6}. $Q\ot\wedge^2Q(j-1),Q\ot S^2Q(j)\in\CC$ for $j=1,2,3,4$.
Consider the exact sequence
$$0\to Q\ot\wedge^2 Q\to V\ot Q^{\ast}\ot Q(1)\to S^2V\ot Q(1)\to Q\ot S^2Q(1)\to 0$$
obtained by tensoring \eqref{S2Q-seq} with $Q(1)$.
Using Steps 4 and 5 we deduce that $Q\ot S^2Q(1)\in\CC$. Note that the subcategory
$\CC$ is admissible, so it is closed under passing to direct summands. Since
$Q\ot S^2Q(1)=S^3Q(1)\oplus S^{(2,1)}Q(1)$, we derive that $S^{(2,1)}Q(1)\in\CC$.
This implies that $Q\ot\wedge^2 Q(1)=Q^{\ast}(2)\oplus S^{(2,1)}Q(1)\in\CC$ (where
$Q^{\ast}(2)\in\CC$ by Step 1). Now we tensor the above exact sequence by
$\OO(1)$ and iterate the above argument.

\vspace*{0.1cm}

\noindent
{\bf Step 7}. $Q\ot S^3Q(2)\in\CC$.
Consider the exact sequence
\begin{equation}\label{S3Q-seq}
0\to Q(-1)\to V\ot\wedge^2 Q(-1)\to S^2V\ot Q^{\ast}\to S^3V\ot\OO\to S^3Q\to 0
\end{equation}
obtained from \eqref{basic-seq}. Tensoring it with $Q(2)$ and using Steps 2, 4 and 6
we deduce the assertion.

\vspace*{0.1cm}

\noindent
{\bf Step 8}. $S^2Q\ot S^2Q(2)\in\CC$. 
We have $S^2Q\ot S^2Q(2)=Q\ot S^3Q(2)\oplus S^{(2,2)}Q(2)$. Hence, by Step 7,
it is enough to check that $S^{(2,2)}Q(2)\in\CC$. But $S^{(2,2)}Q(2)$ is a direct summand in
$\wedge^2Q\ot\wedge^2Q(2)$, so the assertion follows from Step 3.

\vspace*{0.1cm}

\noindent
{\bf Step 9}. $S^4Q(1)\in\CC$.
This follows immediately from the exact sequence
$$0\to\OO\to V\ot Q\to S^2V\ot\wedge^2Q\to S^3V\ot Q^{\ast}(1)\to S^4V\ot\OO(1)\to S^4Q(1)\to 0$$
deduced from \eqref{basic-seq}.

\vspace*{0.1cm}

\noindent
{\bf Step 10}. $\wedge^2Q\ot S^2Q(1)\in\CC$.
Consider the exact sequence
$$0\to\wedge^2 Q^{\ast}\to\wedge^2 V\ot\OO\to V\ot Q\to S^2Q\to 0$$
deduced from \eqref{basic-seq}. Tensoring it with $S^2 Q(2)$ we get the exact sequence
$$0\to \wedge^2 Q\ot S^2Q(1)\to 
\wedge^2 V\ot S^2Q(2)\to V\ot Q\ot S^2Q(2)\to S^2Q\ot S^2Q(2)\to 0.$$
Here all the nonzero terms except for the first one belong to $\CC$ by
Steps 2, 6 and 8, so the assertion follows.

\noindent
{\bf Step 11}. Finally, we are going to deduce that $Q\ot Q\in\CC$.  Tensoring \eqref{S3Q-seq} by
$Q(1)$ we get an exact sequence
$$0\to Q\ot Q\to V\ot Q\ot\wedge^2 Q\to S^2V\ot Q^{\ast}\ot Q(1)\to S^3V\ot Q(1)\to Q\ot S^3Q(1)\to 0.$$
All the nonzero terms except for the first and the last belong to $\CC$ by Steps 4 and 5.
Thus, it is enough to check that $Q\ot S^3Q(1)\in\CC$. We have
$Q\ot S^3Q(1)= S^4Q(1)\oplus S^{(3,1)}Q(1)$. It remains to observe that
$S^4Q(1)\in\CC$ by Step 9, while $S^{(3,1)}Q(1)\in\CC$ as a direct summand of
$\wedge^2 Q\ot S^2Q(1)$ which is in $\CC$ by Step 10.
\ed

\vspace*{0.5cm}

\begin{thm}\label{full-thm}
The exceptional collection on $LG(4,8)$ considered in Theorem \ref{exc-col-thm} is full.
\end{thm}

\Pf . \ Recall that $Q$ is dual to the universal subbundle $\Uu=\Uu_4\sub V\ot\OO$.
Taking the dual of the collection in question we obtain the
collection 
\begin{equation}\label{eq:dualcollection}
(\wedge^2\Uu _4(-4),\Uu _4(-4),\Oo(-4),\ldots,
\wedge^2\Uu _4(-1),\Uu _4(-1),\Oo(-1),\wedge^2\Uu _4,\Uu,E^{\ast},\Oo)
\end{equation}
that generates the admissible triangulated subcategory $\CC^{\ast}\sub\Dd ^{b}(LG(4,8))$. 
It is enough to check that $\CC^{\ast}=\Dd ^{b}(LG(4,8))$. Consider the
diagram with $p$ and $\pi$ being natural projections:
\begin{diagram}
&&{\rm F}_{1,4,8}&&\\
&\ldTo{\pi}&&\rdTo{p}&\\
{\mathbb P}^7&&&& LG(4,8)
\end{diagram}
Here ${\rm F}_{1,4,8}$ is the partial flag variety consisting of pairs
$(l\subset U)$, where $l$ is a line in a Lagrangian subspace $U\sub V$. The
variety ${\rm F}_{1,4,8}$ is naturally embedded into the product
${\mathbb P}^7\times LG(4,8)$. Let us denote by $i:{\rm F}_{1,4,8}\hookrightarrow {\mathbb P}^7\times 
 LG(4,8)$ the natural embedding. Consider the fiber $\pi^{-1}(x)$ over a
point $x$ in $\Pp ^7$. The variety $\pi^{-1}(x)$ is isomorphic to the
Lagrangian Grassmannian $LG(3,6)$. There is a rank
three vector bundle $\Uu _3$ on ${\rm F}_{1,4,8}$ such that its
restriction to any fiber $\pi^{-1}(x)$ is isomorphic to the universal bundle over
this fiber. Recall that the derived category of coherent sheaves on
$\pi^{-1}(x)$ has a full exceptional collection:
\begin{equation}\label{LG36-col}
(\Uu _3|_{\pi^{-1}(x)}\otimes \Oo _{\pi}(-3),\Oo _{\pi^{-1}(x)}(-3),\Uu
_3|_{\pi^{-1}(x)}\otimes \Oo _{\pi}(-2),\Oo _{\pi^{-1}(x)}(-2),\ldots,\Uu _3|_{\pi^{-1}(x)},\Oo).
\end{equation}
Here $\Oo _{\pi}(-1)$ is a line bundle that is isomorphic to ${\rm
  det} \ \Uu _3$. Therefore, by Theorem 3.1 of \cite{Sam},
 the category $\Dd ^{b}({\rm F}_{1,4,8})$ has a semiorthogonal decomposition:
\begin{equation}
\Dd ^{b}({\rm F}_{1,4,8}) = \langle \pi ^{\ast}\Dd ^{b}(\Pp ^7)\otimes
\Uu _3\otimes \Oo _{\pi}(-3),\pi ^{\ast}\Dd ^{b}(\Pp
^7)\otimes  \Oo _{\pi}(-3),\ldots,\pi ^{\ast}\Dd ^{b}(\Pp ^7)\otimes
\Uu _3,\pi ^{\ast}\Dd ^{b}(\Pp ^7)\rangle .
\end{equation}

There is a short exact sequence of vector bundles on ${\rm F}_{1,4,8}$:
\begin{equation}\label{eq:seq1}
0\rightarrow {\pi}^{\ast}\Oo (-1)\rightarrow p^{\ast}\Uu _4\rightarrow
\Uu _3\rightarrow 0.
\end{equation}
Taking determinants we get an isomorphism of line bundles $p^{\ast}\Oo (-1) =
\pi ^{\ast}\Oo (-1)\otimes \Oo _{\pi}(-1)$. 
Therefore, we can replace $\Oo_{\pi}(i)$ by $p^{\ast}\Oo(i)$ in the above semiorthogonal
decomposition. Thus, to prove the statement it
is sufficient to show that all the subcategories 
$$p_{\ast}(\pi ^{\ast}\Dd ^{b}(\Pp ^7)\ot\Uu_3)\ot \Oo(j),
p_{\ast}(\pi ^{\ast}\Dd ^{b}(\Pp ^7))\ot\Oo(j), \ \text{ for }j=0,\ldots,-3$$
belong to $\CC^{\ast}$.

The functor $p_{\ast}\pi^{\ast}:\Dd ^{b}({\mathbb P}^7)\rightarrow \Dd
^{b}(LG(4,8))$ can be computed using the Koszul resolution of
the sheaf $i_{\ast}\Oo _{{\rm F}_{1,4,8}}$ on ${\mathbb P}^7\times LG(4,8)$:
\begin{equation}\label{eq:Koszulres}
0 \rightarrow {\pi}^{\ast}\Oo (-4)\otimes p^{\ast}\Oo (-1) \rightarrow \dots
\rightarrow {\pi}^{\ast}\Oo (-2)\otimes \wedge ^{2}p^{\ast}\Uu _4\rightarrow {\pi}^{\ast}\Oo (-1)\otimes
p^{\ast}\Uu _4\rightarrow \Oo
\rightarrow i_{\ast}{\Oo _{{\rm F}_{1,4,8}}} \rightarrow 0 
\end{equation}

Using this resolution we immediately check the inclusion
$$p_{\ast}(\pi ^{\ast}\Dd ^{b}(\Pp ^7))\ot\Oo(j)\sub 
\lan\Oo(j-1),\wedge ^3\Uu _4(j)=\Uu _4^{\ast}(j-1), \wedge ^2\Uu _4(j),\Uu_4(j),\Oo(j)\ran.$$
By Lemma \ref{full-lem}(i), for $j=-3,\ldots,0$ the right-hand side belongs to $\CC^{\ast}$.

Next, using the sequence \eqref{eq:seq1} we see that to prove the inclusion
$p_{\ast}(\pi ^{\ast}\Dd ^{b}(\Pp ^7)\ot\Uu_3)\ot \Oo(j)\sub\CC^{\ast}$ it is enough to check
that
$$\lan\Uu_4(j-1), \Uu_4\ot\Uu _4^{\ast}(j-1), \Uu_4\ot\wedge ^2\Uu _4(j),\Uu_4\ot\Uu_4(j),\Uu_4(j)\ran
\sub\CC^{\ast}$$
for $j=-3,\ldots,0$. It remains to apply Lemma \ref{full-lem} (and dualize).
\ed

\noindent{\it Another version of the proof}.
We can simplify computations in the above argument by using a different semiorthogonal
decomposition of $\Dd^b({\rm F}_{1,4,8}$):
$$\Dd ^{b}({\rm F}_{1,4,8}) = \langle \pi ^{\ast}\Dd ^{b}(\Pp ^7)\otimes
\Uu _3\otimes p^{\ast}\Oo(-3),\pi ^{\ast}\Dd ^{b}(\Pp
^7)\otimes  p^{\ast}\Oo(-3),\ldots,
\pi ^{\ast}\Dd ^{b}(\Pp ^7),\pi ^{\ast}\Dd ^{b}(\Pp ^7)\otimes
\Uu _3^{\ast}\rangle.$$
The restriction of this decomposition to the fiber $\pi^{-1}(x)\simeq LG(3,6)$ is  
the exceptional collection obtained
from collection \eqref{LG36-col} by the right mutation of $\Uu_3|_{\pi^{-1}(x)}$ through
$\OO$. In the same way as above we check that
$$p_{\ast}(\pi ^{\ast}\Dd ^{b}(\Pp ^7))\ot\Oo(j)\sub\CC^{\ast} \ \text{ for }j=-3,\ldots,0,$$
$$p_{\ast}(\pi ^{\ast}\Dd ^{b}(\Pp ^7)\ot\Uu_3)\ot \Oo(j)\sub\CC^{\ast} \ \text{ for }j=-1,-2,-3, \ \text{and}$$
$$p_{\ast}(\pi^{\ast}\Oo(i)\ot\Uu_3^{\ast})\in\CC^{\ast} \ \text{ for }i=-6,\ldots,0.$$
The point is that this will only require using (easy) Steps 1,2,4,5 and 6 of Lemma \ref{full-lem}.
Thus, if we consider the semiorthogonal decomposition
$$\Dd^b({\rm F_{1,4,8}})=\lan \AA, \lan\pi^{\ast}\Oo(1)\ot\Uu_3^{\ast}\ran\ran,$$
where $\AA=\lan\pi^{\ast}\Oo(1)\ot\Uu_3^{\ast}\ran^{\perp}$, then
$p_{\ast}\AA\in\CC^{\ast}$. By adjointness, it follows that for an object 
$E\in\Dd^b(LG(4,8))$ such that $\Hom(E,\CC^{\ast})=0$, 
one has $p^{\ast}E\in\lan\pi^{\ast}\Oo(1)\ot\Uu_3^{\ast}\ran$,
i.e., $p^{\ast}E\simeq V^{\bullet}\ot \pi^{\ast}\Oo(1)\ot\Uu_3^{\ast}$ for a graded vector space
$V^{\bullet}$. Hence, $E\simeq V^{\bullet}\ot p_{\ast}(\pi^{\ast}\Oo(1)\ot\Uu_3^{\ast})$.
Finally, using resolution \eqref{eq:Koszulres} and the dual of sequence \eqref{eq:seq1}
one can compute that 
$$p_{\ast}(\pi^{\ast}\Oo(1)\ot\Uu_3^{\ast})\simeq\we^2\Uu_4^{\ast}\simeq\we^2\Uu_4(1).$$
Thus, $E\simeq V^{\bullet}\ot\we^2\Uu_4(1)$. But $\Hom^*(\we^2\Uu_4(1),\we^2\Uu_4(-4))\neq 0$
by Serre duality, so the condition $\Hom^*(E,\CC^{\ast})=0$ implies that $V^{\bullet}=0$.
\ed

\section{The case of $LG(5,10)$}

In this section we assume that $n=5$ (so $V$ is $10$-dimensional).
It turns out that in this case the exceptional bundles constructed so far
do not generate the entire derived category $\Dd^b(LG(5,10))$. We are going to construct
another exceptional bundle on $LG(5,10)$ starting from the bundle $T=S^{(3,1,1)}Q$.
Let us denote by $\om_i$ the $i$th fundamental weight of the root system $C_5$.
For a dominant weight $\la$ we denote by $V(\la)$ the corresponding irreducible representation of
$\Sp(10)$ (for example, $V(\om_1)=V$, $V(\om_2)=\we^2 V/k$, $V(2\om_1)=S^2V$).

\begin{lem}\label{T-lem} 
Assume that $n=5$. 

\noindent
(i) $\Hom^*(\we^iQ,T(j))=0$ for $i\in [0,3]$, $j\in [-5,-1]$. Also, $\Hom^*(T,\OO)=0$.

\noindent
(ii) $\Hom^*(R_1,T(j))=0$ for $j\in [-5,-1]$.

\noindent
(iii) $\Hom^*(T,T(-3))=0$.

\noindent
(iv) $\Hom^i(T,T)=0$ for $i>2$, $\Hom^2(T,T)=V(2\om_1+\om_2)\oplus V(\om_1+\om_3)$,
$\Hom^1(T,T)=V^{\ot 2}/k\oplus S^2V$, $\Hom^0(T,T)=k$.

\noindent
(v) $\Hom^i(\we^3Q,T)=0$ for $i>0$ and $\Hom^0(\we^3Q,T)=S^2V$. Also,
$\Hom^i(T,\we^3Q)=0$ for $i\neq 1,2$, $\Hom^1(T,\we^3Q)=k$ and $\Hom^2(T,\we^3Q)=\we^2V/k$.

\noindent
(vi) $\Hom^i(\we^2Q,T)=0$ for $i>0$ and $\Hom^0(\we^2Q,T)=V(3\om_1)\oplus V(\om_1+\om_2)$.
Also, $\Hom^i(T,\we^2Q)=0$ for $i\neq 2$ and $\Hom^2(T,\we^2Q)=V$.

\noindent
(vii) $\Hom^i(Q,T)=0$ for $i>0$ and $\Hom^0(Q,T)=V(2\om_1+\om_2)\oplus V(\om_1+\om_3)$.
Also, $\Hom^i(T,Q)=0$ for $i\neq 2$ and $\Hom^2(T,Q)=k$.

\noindent
(viii) $\Hom^i(R_1,T)=0$ for $i>1$, $\Hom^1(R_1,T)=V(2\om_1+\om_2)\oplus V(\om_1+\om_3)$ and 
$\Hom^0(R_1,T)=V^{\ot 2}/k$. Also, $\Hom^i(T,R_1)=0$ for $i\neq 1,2$, $\Hom^1(T,R_1)=k$ and
$\Hom^2(T,R_1)=V^{\ot 2}/k$.

\noindent
(ix) $\Hom^*(T,S^2Q^{\ast})=0$.

\noindent
(x) $\Hom^i(T,S^3Q^{\ast})=0$ for $i\neq 4$.

\noindent
(xi) $\Hom^i(T,R_3)=0$ for $i\neq 1$.

\noindent
(xii) $\Hom^i(T,\we^3Q\ot Q^{\ast})=0$ for $i\neq 2$.

\noindent
(xiii) $\Hom^4(T,\we^3Q\ot\we^2Q^{\ast})=0$.
\end{lem}

The proof is a straightforward application of the Bott's theorem. 
By part (viii) of the above Lemma, we have a canonical nonsplit extension of vector bundles
\begin{equation}\label{PT-seq}
0\to R_1\to P\to T\to 0
\end{equation}

\begin{lem}\label{P-lem} 
(i) The map $\Hom^1(R_1,Q)\to\Hom^2(T,Q)$ induced by \eqref{PT-seq} is an isomorphism.

\noindent
(ii) The map $\Hom^1(R_1,\we^2Q)\to\Hom^2(T,\we^2Q)$ induced by \eqref{PT-seq} is an isomorphism.

\noindent
(iii) The map $\Hom^1(R_1,R_1)\to\Hom^2(T,R_1)$ induced by \eqref{PT-seq} is an isomorphism.

\noindent
(iv) The map $\Hom^1(R_1,T)\to\Hom^2(T,T)$ induced by \eqref{PT-seq} is an isomorphism, while
the map $\Hom^0(R_1,T)\to\Hom^1(T,T)$ is injective.

\noindent
(v) One has $\Hom^*(P,Q)=\Hom^*(P,\we^2Q)=\Hom^*(P,R_1)=\Hom^{>1}(P,P)=0$
and $\Hom^1(P,P)=S^2V$, $\Hom^0(P,P)=k$. 
Also, $\Hom^i(P,\we^3Q)=0$ for $i\neq 1$ and $\Hom^1(P,\we^3 Q)=k$.
\end{lem}

\Pf . (i) We have to check that the natural map
$$\Hom^1(R_1,Q)\ot\Hom^1(T,R_1)\to\Hom^2(T,Q)$$
is an isomorphism. Note that both sides are one-dimensional (see Lemma \ref{Bott-lem}(iii) and
Lemma \ref{T-lem}(vii),(viii)), so it is enough to
check that this map is nonzero. We have natural embeddings
$S^2Q^{\ast}\to R_1^{\ast}\ot Q$ and $S^2Q^{\ast}\to T^{\ast}\ot R_1$ inducing
isomorphisms on $H^1$. Let us consider the induced map
$$\a:S^2Q^{\ast}\ot S^2Q^{\ast}\to T^{\ast}\ot Q.$$
Note that 
$$S^2Q^{\ast}\ot S^2Q^{\ast}=S^4Q^{\ast}\oplus S^{(2,2)}Q^{\ast}\oplus S^{(3,1)}Q^{\ast},$$
where the first two terms have zero cohomology while the last term has one-dimensional $H^2$.
Thus, it is enough to check that the restriction of $\a$ to $S^{(3,1)}Q^{\ast}$ is nonzero and
that the natural map
$$H^1(S^2Q^{\ast})\ot H^1(S^2Q^{\ast})\to H^2(S^2Q^{\ast}\ot S^2Q^{\ast})$$
between one-dimensional spaces is nonzero. Let us start by splitting the exact sequence 
\eqref{S2Q*-seq} into two short exact sequences
\begin{equation}\label{K-eq}
0\to S^2Q^{\ast}\to S^2V\ot\OO\to K\to  0
\end{equation}
\begin{equation}\label{K2-eq}
0\to K\to V\ot Q\to \we^2 Q\to 0
\end{equation}
Then \eqref{K-eq} induces the surjections
$H^0(K)\to H^1(S^2Q^{\ast})$ and 
$H^1(K\ot S^2Q^{\ast})\to H^2(S^2Q^{\ast}\ot S^2Q^{\ast})$
(by the vanishing of $H^1(\OO)$ and $H^2(S^2Q^{\ast})$).
Hence, it is enough to check that the natural map
$$H^0(K)\ot H^1(S^2Q^{\ast})\to H^1(K\ot S^2Q^{\ast})$$
is an isomorphism (note that both sides are isomorphic to $S^2V\oplus k$).
Now the sequence \eqref{K2-eq} induces embeddings
$H^0(K)\to V\ot H^0(Q)$ and $H^1(K\ot S^2Q^{\ast})\to V\ot H^1(Q\ot S^2Q^{\ast})$
(by the vanishing of $H^0(\we^2Q\ot S^2Q^{\ast})$). Hence, we are reduce to proving
that the map 
$$H^0(Q)\ot H^1(S^2Q^{\ast})\to H^1(Q\ot S^2Q^{\ast})$$
is an isomorphism. But this follows from the exact sequence \eqref{basic-seq} and the vanishing
of $H^*(Q^{\ast}\ot S^2Q^{\ast})$.

It remains to check that the restriction of $\a$ to $S^{3,1}Q^{\ast}\sub S^2Q^{\ast}\ot S^2Q^{\ast}$
is nonzero (where we can just think of $Q$ as a vector space). 
Let us view $T$ (resp., $R_1$) as the image of the Koszul differential
$S^2Q\ot\we^3Q\to S^3Q\ot \we^2Q$ (resp., $Q\ot\we^2Q\to S^2Q\ot Q$). 
Then the embedding $S^2Q^{\ast}\hra T^{\ast}\ot R_1$ corresponds to the composed map
\begin{equation}\label{TR1-eq}
S^2Q^{\ast}\ot T\to S^2Q^{\ast}\ot S^3Q\ot \we^2Q\to Q\ot\we^2Q\to R_1,
\end{equation}
where the second arrow is induced by the natural map $S^2Q^{\ast}\ot S^3Q\to Q$.
On the other hand, the embedding $S^2Q^{\ast}\hra R_1^{\ast}\ot Q$ corresponds to
the natural map $S^2Q^{\ast}\ot R_1\to Q$ induced by the embedding $R_1\to S^2Q\ot Q$.
Thus, $\a$ corresponds to the composed map
$$
\a':S^2Q^{\ast}\ot S^2Q^{\ast}\ot T\to S^2Q^{\ast}\ot S^2Q^{\ast}\ot S^3Q\ot\we^2Q\to
S^2Q^{\ast}\ot Q\ot\we^2Q\to S^2Q^{\ast}\ot S^2Q\ot Q\to Q,
$$
where the third arrow is induced by the Koszul differential.
Let us choose a basis $e_1,\ldots,e_n$ for $Q$
and define an element $t\in T$ by
$$t=e_4^2e_1\ot (e_2\we e_3)+e_4^2e_2\ot (e_3\we e_1)+e_4^2e_3\ot (e_1\we e_2),$$
where we view $T$ as a subbundle in $S^3Q\ot \we^2Q$. Then one can compute 
the induced functional on $S^2Q^{\ast}\ot S^2Q^{\ast}$
$$x\mapsto \lan \a'(x\ot t),e_3\ran=2\lan x, (e_0e_2)\we(e_0e_1)\ran,$$
where $x\in S^2Q^{\ast}\ot S^2Q^{\ast}$. Now we observe that $S^{(3,1)}Q$ can be identified
with the image of $\we^2(S^2Q)$ under 
the natural map $\b:S^2Q\ot S^2Q\to S^3Q\ot Q$ given by 
$$\b(f\ot(v_1v_2))=(fv_1)\ot v_2+(fv_2)\ot v_1.$$
Finally we compute that 
$$\b((e_0e_2)\we(e_0e_1))=(e_4^2e_2)\ot e_1-(e_4^2e_1)\ot e_2\neq 0,$$
which finishes the proof.

\noindent
(ii) Since both the source and the target are isomorphic to $V$, it is enough to check surjectivity.
Furthermore, it suffices to prove that the composition map
$$\Hom^1(T,R_1)\ot\Hom^1(R_1,Q)\ot\Hom^0(Q,\we^2Q)\to\Hom^2(T,\we^2Q)$$
is surjective. By part (i), this reduces to surjectivity of the composition map
$$\Hom^2(T,Q)\ot\Hom^0(Q,\we^2Q)\to\Hom^2(T,\we^2Q).$$
Looking at the exact sequence \eqref{K2-eq}, we see that this would follow from
the vanishing of $\Hom^3(T,K)$. But this vanishing follows from the exact sequence \eqref{K-eq} 
since $\Hom^*(T,\OO)=\Hom^*(T,S^2Q^{\ast})=0$ (see Lemma \ref{T-lem}(i),(ix)).

\noindent
(iii) Both the source and the target are isomorphic to $V^{\ot 2}/k$ (see Lemma 
\ref{Bott-lem}(iv) and Lemma \ref{T-lem}(viii)), so it suffices to
check surjectivity. By part (ii), it is enough to prove that the map
\begin{equation}\label{TR1-map}
\Hom^2(T,\we^2Q)\ot V\to\Hom^2(T,R_1)
\end{equation}
is surjective. Let us first check that $S^2V\sub\Hom^2(T,R_1)$ is in the image.
The exact sequence \eqref{basic-seq} induces a long exact sequence
$$\ldots\to \Hom^2(T,\we^2Q\ot Q^{\ast})\to \Hom^2(T,\we^2Q)\ot V\stackrel{f}{\to}
\Hom^2(T,\we^2Q\ot Q)\to\ldots$$
Using the Bott's theorem one can check that $\Hom^2(T,\we^2Q\ot Q^{\ast})$ does not contain
any factors isomorphic to $S^2V$, so the restriction of $f$ to $S^2V$ is an embedding.
On the other hand, 
$$\Hom^2(T,\we^2Q\ot Q)=\Hom^2(T,R_1)\oplus\Hom^2(T,\we^3Q),$$
where the second factor is $\we^2V/k$, hence, $S^2V$ projects nontrivially to $\Hom^2(T,R_1)$.
It remains to check that $\we^2V/k\sub\Hom^2(T,R_1)$ is in the image of the map \eqref{TR1-map}.
It suffices to prove that it is in the image of the map
$$\Hom^2(T,Q\ot Q)\ot V\to\Hom^2(T,R_1),$$
or even
$$\Hom^2(T,Q)\ot H^0(\we^2Q)\to\Hom^2(T,R_1).$$
We have a natural map
$$\ga:S^2Q^{\ast}\ot \we^2Q^{\ast}\to T^{\ast}\ot Q,$$
such that its composition with the embedding $T^{\ast}\ot Q\hra S^2Q^{\ast}\ot \we^3Q^{\ast} \ot Q$
(induced by the surjection $S^2Q\ot\we^3Q\to T$) is the identity map on $S^2Q^{\ast}$ tensored
with the natural embedding $\we^2Q^{\ast}\to\we^3Q^{\ast}\ot Q$. Note that this implies that
$\ga$ itself is an embedding. Hence, $\ga$ induces an isomorphism on $H^2$.
Next, we claim that the composition map
$$H^2(S^2Q^{\ast}\ot\we^2Q^{\ast})\ot H^0(\we^2Q)\to H^2(S^2Q^{\ast}\ot\we^2Q^{\ast}\ot\we^2Q)$$
is surjective. Indeed, it is enough to check this with $H^0(\we^2Q)$ replaced by $\we^2 V$.
Then the exact sequence 
$$0\to S^2Q^{\ast}\to V\ot Q^{\ast}\to \we^2 V\ot \OO\to\we^2Q\to 0$$
shows that this follows from the vanishing of
$H^3(S^2Q^{\ast}\ot\we^2Q^{\ast}\ot Q^{\ast})$ and $H^4(S^2Q^{\ast}\ot\we^2Q^{\ast}\ot S^2Q^{\ast})$,
both of which are easily checked using the Bott's theorem.
Now it remains to prove that the composed map
$$S^2Q^{\ast}\ot\we^2Q^{\ast}\ot\we^2Q\stackrel{\ga\ot\id}{\to} T^{\ast}\ot Q\ot\we^2Q\to T^{\ast}\ot R_1$$
induces an embedding on $H^2$. It is enough to prove that the kernel of this map is $S^2Q^{\ast}$.
Using the embedding of $T^{\ast}$ into $S^2Q^{\ast}\ot \we^3Q^{\ast}$ this reduces to checking
that the composition of the natural maps
$$\we^2Q^{\ast}\ot\we^2Q\to \we^3Q^{\ast}\ot Q\ot\we^2Q\to\we^3Q^{\ast}\ot R_1$$ 
has $\OO$ as a kernel. Replacing this map by its composition with the
embedding $\we^3Q^{\ast}\ot R_1\hra\we^3Q^{\ast}\ot S^2Q\ot Q$ we see that it is enough
to prove the following fact from linear algebra.
Suppose we have a linear map $A:\we^2Q\to \we^2Q$ such that the induced map
$$\we^3Q\to Q\ot\we^2Q\stackrel{\id\ot A}{\to} Q\ot\we^2Q\stackrel{d}{\to} S^2Q\ot Q$$
is zero, where $d$ is Koszul differential. Then $A$ is proportional to identity. To prove this statement 
we recall
that the kernel of $d$ is exactly $\we^3Q\sub Q\ot\we^2Q$. Thus, the condition on $A$ is that
$\id_Q\ot A$ preserves $\we^3Q\sub Q\ot\we^2Q$. Let us fix some basis $(e_i)$ of $Q$ and
let $\partial_i:\we^3Q\to\we^2Q$ be the odd partial derivatives corresponding to the dual basis
of $Q^{\ast}$. Consider
$$e_1\ot A(e_2\we e_3)+e_2\ot A(e_3\we e_1)+ e_3\ot A(e_1\we e_2)=\eta\in\we^3Q\sub Q\ot \we^2Q.$$
Contracting with $e_3^*$ in the first factor of the tensor product $Q\ot\we^2$ we obtain
$A(e_1\we e_2)=\partial_3\eta$. Hence, 
$\partial_3A(e_1\we e_2)=\partial_3^2\eta=0$. In a similar way
$\partial_i A(e_1\we e_2)=0$ for $i>2$. It follows that $A(e_1\we e_2)$ is proportional to $e_1\we e_2$.
Thus, for every pair of elements $x,y\in Q$, $A(x\we y)$ is proportional to $x\we y$.
This implies that $A$ is proportional to identity.


\noindent
(iv) We have $\Hom^1(R_1,T)\simeq\Hom^2(T,T)$ (see Lemma \ref{T-lem}(iv),(viii)), so for the first assertion
it is enough to check the surjectivity.
By part (ii), it suffices to check that the map
$$\Hom^2(T,\we^2Q)\ot\Hom^0(\we^2Q,T)\to\Hom^2(T,T)$$
is surjective. Furthermore, it is enough to prove that the map
$$\Hom^2(T,\we^2Q)\ot\Hom^0(\we^2Q,\we^3Q)\ot\Hom^0(\we^3Q,T)\to\Hom^2(T,T)$$
is surjective. We are going to do this in two steps: first, we'll check that the map
\begin{equation}\label{Twe23-eq}
\Hom^2(T,\we^2Q)\ot V\to \Hom^2(T,\we^3Q)
\end{equation}
is surjective, and then we will show the surjectivity of
\begin{equation}\label{Twe2T-eq}
\Hom^2(T,\we^3Q)\ot\Hom^0(\we^3Q,T)\to\Hom^2(T,T)
\end{equation}
From the exact sequence \eqref{basic-seq} we get the following long exact sequence
$$0\to S^3Q^{\ast}\to S^3V\ot\OO\to S^2V\ot Q\to V\ot\we^2Q\to\we^3Q\to 0.$$
Thus, the surjectivity of \eqref{Twe23-eq} follows from the vanishing of
$\Hom^3(T, Q)$, $\Hom^4(T,\OO)$ and $\Hom^5(T,S^3Q^{\ast})$ (see Lemma \ref{T-lem}(i),(vii),(x)).
To deal with \eqref{Twe2T-eq}
we use the natural embedding $S^2Q\to\we^3Q^{\ast}\ot T$ inducing an isomorphism
on $H^0$. Note also that since $\we^3Q\ot S^2Q\simeq T\oplus R_3$, Lemma \ref{T-lem}(xi) implies
that the projection $T^{\ast}\ot\we^3Q\ot S^2Q\to T^{\ast}\ot T$ induces an isomorphism on $H^2$.
Thus, we are reduced to showing the surjectivity of
$$\Hom^2(T,\we^3 Q)\ot H^0(S^2Q)\to \Hom^2(T,\we^3Q\ot S^2Q).$$
It suffices to prove the surjectivity of the maps
$$\Hom^2(T,\we^3 Q)\ot V\to\Hom^2(T,\we^3 Q\ot Q),$$
$$\Hom^2(T,\we^3 Q\ot Q)\ot V\to\Hom^2(T,\we^3 Q\ot S^2Q).$$
The exact sequence \eqref{basic-seq} shows that the surjectivity of the first map follows from
the vanishing of $\Hom^3(T,\we^3 Q\ot Q^{\ast})$ (see Lemma \ref{T-lem}(xii)). Similarly, for
the second map we use the exact sequence 
$$0\to\wedge^2 Q^{\ast}\to\wedge^2 V\ot\OO\to V\ot Q\to S^2Q\to 0$$
 along with the vanishing of $\Hom^3(T,\we^3Q)$
and $\Hom^4(T,\we^3Q\ot\we^2Q^{\ast})$ (see Lemma \ref{T-lem}(v),(xiii)).

Now let us prove the injectivity of the map $\Hom^0(R_1,T)\to\Hom^1(T,T)$. We have a natural
embedding $S^2Q\to R_1^{\ast}\ot T$ inducing isomorphism on $H^0$ and an embedding
$S^2Q^{\ast}\to T^{\ast}\ot R_1$ inducing isomorphism on $H^1$. We claim that the
composed map 
\begin{equation}\label{S2T-map}
S^2Q\ot S^2Q^{\ast}\to T^{\ast}\ot T
\end{equation} 
induces an embedding on $S^{(2,0,0,0,-2)}Q\sub S^2Q\ot S^2Q^{\ast}$. 
To prove this we can replace $Q$ by a vector space with a basis $e_1,\ldots,e_5$. Let
$e_1^*,\ldots,e_5^*$ be the dual basis of $Q^{\ast}$. 
It is enough to check that the lowest weight vector $e_1^2\ot (e_5^*)^2$ maps to a nonzero element
of $T^{\ast}\ot T$ under \eqref{S2T-map}. By definition, this endomorphism of $T$ is the composition
of the map
$$T\to S^3Q\ot\we^2Q\stackrel{\partial_5^2\ot\id}{\to} Q\ot\we^2Q\to R_1$$
with the map
$$R_1\to Q\ot\we^2Q\stackrel{e_1^2}{\to} S^2Q\ot Q\ot\we^2Q\to S^3Q\ot\we^2Q\to T.$$
Viewing $T$ as a direct summand of $S^2Q\ot\we^3Q$ we obtain from the first (resp., second) map
a map $f:S^2Q\ot\we^3Q\to R_1$ (resp., $g:R_1\to S^2Q\ot\we^3Q$).
Identifying $R_1$ with $Q\ot\we^2Q/\we^3Q$ we can write
$$f(t\ot (x\we y\we z))=\partial_5^2(tx)\ot(y\we z)+\partial_5^2(ty)\ot(z\we x)+\partial_5^2(tz)\ot(x\we y)
\mod\we^3Q,$$
$$g(x\ot (y\we z)\mod\we^3Q)=2(e_1x)\ot(e_1\we y\we z)+(e_1y)\ot(e_1\we x\we z)+
(e_1z)\ot(e_1\we y\we x),$$
where $t\in S^2Q$ and $x,y,z\in Q$ (for appropriate rescaling of $g$). Hence,
$$gf((e_4e_5)\ot(e_1\we e_2\we e_5))=2g(e_4\ot(e_1\we e_2)\mod\we^3Q)=2e_1^2\ot(e_1\we e_4
\we e_2)\neq 0.$$
Thus, the map \eqref{S2T-map} induces
an embedding on $H^1$. So we are reduced to checking that the natural map
$$H^0(S^2Q)\ot H^1(S^2Q^{\ast})\to H^1(S^2Q\ot S^2Q^{\ast})$$
is an isomorphism. Since both sides are isomorphic to $S^2V$, it is enough to prove surjectivity.
The exact sequence \eqref{S2Q-seq}
 shows that this follows from the vanishing of $H^2(Q^{\ast}\ot S^2Q^{\ast})$
and $H^3(\we^2Q^{\ast}\ot S^2Q^{\ast})$, which can be checked using the Bott's theorem.

\noindent 
(v) The vanishing of $\Hom^*(P,Q)$, $\Hom^*(P,\we^2Q)$, $\Hom^*(P,R_1)$ 
follow from directly from parts (i)-(iv) along with the computation of the relevant spaces in
Lemmas \ref{Bott-lem} and \ref{T-lem}. Similarly, we derive that $\Hom^0(P,T)=k$,
$\Hom^1(P,T)=S^2V$ and $\Hom^i(P,T)=0$ for $i>1$. Now one computes $\Hom^*(P,P)$
by applying the functor $\Hom(P,?)$
to the exact sequence \eqref{PT-seq} and using the vanishing of $\Hom^*(P,R_1)$.
To compute $\Hom^*(P,\we^3Q)$ it remains to check that the map
$$\Hom^1(R_1,\we^3Q)\to\Hom^2(T,\we^3Q)$$
induced by \eqref{PT-seq} is an isomorphism. Since both sides are isomorphic to
$\we^2V/k$, it is enough to prove surjectivity.
But this follows immediately from part (ii) along
with the surjectivity of the map \eqref{Twe23-eq} proved in part (iv).
\ed

By part (v) of the above Lemma, we have a canonical nonsplit extension of vector bundles
\begin{equation}\label{PG-seq}
0\to\we^3 Q\to G\to P\to 0
\end{equation}

\begin{thm}\label{G-thm} 
The vector bundle $G$ is exceptional and $\Hom^*(G,\we^3Q)=0$.
\end{thm}

\Pf . First, applying the functor $\Hom(?,\we^3Q)$ to the sequence \eqref{PG-seq}
and using Lemma \ref{P-lem}(v)
we find that $\Hom^*(G,\we^3 Q)=0$. Next, applying the functor $\Hom(G,?)$ to this
sequence we derive isomorphisms $\Hom^i(G,G)\simeq\Hom^i(G,P)$.
Recall that $\Hom^*(\we^3Q,R_1)=0$ by Lemma \ref{R-van-lem}. Hence,
applying the functor $\Hom(\we^3Q,?)$ to the sequence \eqref{PT-seq} and using
Lemma \ref{T-lem}(v) we obtain that $\Hom^i(\we^3Q,P)=0$ for $i>0$ and $\Hom^0(\we^3Q,P)=S^2V$.
Thus, using the sequence \eqref{PG-seq} again along with the computation of $\Hom^*(P,P)$
(see Lemma \ref{P-lem}(v))
we see that it is enough to check that the natural map
$$\Hom^0(\we^3Q,P)\ot\Hom^1(P,\we^3Q)\to\Hom^1(P,P)$$
is an isomorphism. Since $\Hom^*(P,R_1)=\Hom^*(\we^3Q,R_1)=0$ (see Lemma \ref{P-lem}(v)),
the exact sequence \eqref{PT-seq} gives an isomorphism of the above map with
$$\Hom^0(\we^3Q,T)\to\Hom^1(P,T)$$
induced by a nonzero element in $\Hom^1(P,\we^3Q)$.
Since the natural map $\Hom^1(T,\we^3Q)\to\Hom^1(P,\we^3Q)$ is an isomorphism (as we have
seen in the proof of Lemma \ref{P-lem}(v)), the above map factors
as the composition of the map
$$\Hom^0(\we^3Q,T)\stackrel{f}{\to}\Hom^1(T,T)$$
induced by a nonzero element in $\Hom^1(T,\we^3Q)$ followed by the map $h$ in the exact
sequence
$$0\to \Hom^0(R_1,T)\stackrel{g}{\to} \Hom^1(T,T)\stackrel{h}{\to} \Hom^1(P,T)\to 0.$$
Thus, it is enough to check that the images of the maps $f$ and $g$
are complementary in $\Hom^1(T,T)$. Since $\Hom^0(R_1,T)=V^{\ot 2}/k$,
$\Hom^0(\we^3Q,T)=S^2V$, while $\Hom^1(T,T)=V^{\ot 2}/k\oplus S^2V$ (see Lemma
\ref{T-lem}(iv),(v),(viii)), it suffices to prove that the images of $S^2V$ under $f$ and $g$
have trivial intersection. Note that we have a natural embedding $S^2Q\to R_1^{\ast}\ot T$ (resp.,
$S^2Q\to \we^3Q^{\ast}\ot T$)
inducing an embedding of $S^2V$ into $\Hom^0(R_1,T)$ (resp., into $\Hom^0(\we^3Q,T)$).
On the other hand, a nonzero element in $\Hom^1(T,\we^3Q)$ is the image of the nonzero
element in $H^1(S^2Q^{\ast})$ with respect to the embedding $S^2Q^{\ast}\to T^{\ast}\ot\we^3Q$.
Furthermore, we have seen in the end of the proof of Lemma \ref{P-lem}(iv) that the natural
map $H^0(S^2Q)\ot H^1(S^2Q^{\ast})\to H^1(S^2Q\ot S^2Q^{\ast})$ is an isomorphism.
Thus, it is enough to prove that the natural maps
$$\a:S^2Q^{\ast}\ot S^2Q\to (T^{\ast}\ot\we^3Q)\ot (\we^3Q^{\ast}\ot T)\to T^{\ast}\ot T \text{  and}$$
$$\b:S^2Q^{\ast}\ot S^2Q\to (T^{\ast}\ot R_1)\ot (R_1^{\ast}\ot T)\to T^{\ast}\ot T$$
induce linear independent maps on $H^1$.
In fact, since $H^1(S^2Q^{\ast}\ot S^2Q)$ comes from the summand 
$S^{(2,0,0,0,-2)}Q\sub S^2Q^{\ast}\ot S^2Q$, generated by the lowest weight vector 
$v=(e_5^*)^2\ot e_1^2$ (where $(e_i)$ is the basis of $Q$),
it suffices to check that $\a(v)$ and $\b(v)$ are not proportional in $T^{\ast}\ot T$.
Recall that in the proof of Lemma \ref{T-lem}(iv) we have 
constructed the maps $f:S^2Q\ot\we^3Q\to R_1$ and $g:R_1\to S^2Q\ot\we^3Q$
such that $gf$ is a multiple of the composition
$$S^2Q\ot\we^3Q\to T\stackrel{\b(v)}{\to} T\to S^2Q\ot\we^3Q.$$
On the other hand, $\a(v)$ is given by the following composition
$$T\to S^2Q\ot\we^3Q\stackrel{\partial_5^2}{\to}\we^3Q\stackrel{e_1^2}{\to} S^2Q\ot\we^3Q\to T.$$
Let us denote by $\pi:S^2Q\ot\we^3Q\to S^2Q\ot\we^3Q$ the projection with the image $T$,
given by 
$$\pi(ab\ot(x\we y\we z))=\frac{3}{5}ab\ot(x\we y\we z)+
\left(ax\ot (b\we y\we z)+bx\ot (a\we y\we z)+c.p.(x,y,z)\right),$$
where $a,b,x,y,z,\in Q$, the omitted terms $c.p.(x,y,z)$ are obtained by cyclically permuting $x,y,z$.
Then we are reduced to checking that $gf$ is not proportional to the composition
$$h:S^2Q\ot\we^3Q\stackrel{\pi}{\to}S^2Q\ot\we^3Q\stackrel{\partial_5^2}{\to}\we^3Q\stackrel{e_1^2}{\to} S^2Q\ot\we^3Q\stackrel{\pi}{\to}S^2Q\ot\we^3Q.$$
To this end we compute
$$\frac{1}{2}gf(e_4e_5\ot (e_2\we e_3\we e_5)=f(e_4\ot(e_2\we e_3))=
2e_1e_4\ot(e_1\we e_2\we e_3)-e_1e_2\ot (e_1\we e_3\we e_4)+e_1e_3\ot (e_1\we e_2\we e_4),$$
$$\frac{25}{2}h(e_4e_5\ot (e_2\we e_3\we e_5)=
3e_1^2\ot(e_2\we e_3\we e_4)+2e_1e_4\ot(e_1\we e_2\we e_3)+2e_1e_2\ot(e_1\we e_3\we e_4)
-2e_1e_3\ot (e_1\we e_2\we e_4),$$
which are clearly not proportional.
\ed

\begin{lem}\label{end-R-lem} 
On $LG(5,10)$ one has $\Hom^*(R_1,R_1(i))=0$ for $i\in [-5,-1]$.
\end{lem}

\Pf . The proof is similar to that of Lemma \ref{Bott-lem}(iv) and is left to the reader.
\ed

\begin{thm}\label{L5-col-thm} 
Let us consider the following two blocks:
$$\AA=(\OO,Q,\we^2 Q,F_1,\we^3 Q,G) \ \text{ and }\ \ \BB=(\OO,Q,\we^2 Q,F_1,\we^3 Q).$$
Then $(\AA,\BB(1),\BB(2),\AA(3),\BB(4),\BB(5))$ is a full exceptional collection in $\Dd^b(LG(5,10))$.
\end{thm}

\Pf . The required semiorthogonality conditions not involving $G$ follow from the fact that $F_1$
is the right mutation of $E_1$ through $\we^2Q$ and from Lemmas \ref{gen-exc-lem}, \ref{split-lem},
\ref{E-van-lem} and \ref{end-R-lem}. 
Using Serre duality and sequences \eqref{PT-seq} and \eqref{PG-seq} we can
reduce all the remaining semiorthogonality conditions to Lemmas
\ref{T-lem}, \ref{P-lem} and Theorem \ref{G-thm} (for $\Hom^*(G(3),G)=0$ we need in addition
the vanishing of $\Hom^*(\we^3Q(3),\we^3Q)$ and $\Hom^*(R_1(3),R_1)$ that follows from
Lemmas \ref{gen-exc-lem} and \ref{end-R-lem}). 

Now let us prove that our exceptional collection is full. 
Following the method of proof of Theorem \ref{full-thm} (involving the partial isotropic flag manifold
${\rm F}_{1,5,10}$ and the relative analog of our collection for $LG(4,8)$)
one can reduce this to checking that the subcategory $\CC$
generated by our exceptional collection contains the subcategories
$$\PP\ot\Oo(j), \PP\ot Q(j), \PP\ot \we^2 Q(j), \PP\ot Q\ot Q, \PP\ot Q\ot\we^2 Q,$$
where $j=0,\ldots,4$ and $\PP=\lan \OO, Q,\we^2 Q,\we^3 Q, Q^{\ast}(1), \OO(1)\ran$.
This gives the following list of objects that have to be in $\CC$:

\noindent
(i) $\OO(j)$, $Q(j)$, $\we^2Q(j)$ for $j=0,\ldots,5$;

\noindent
(ii) $Q\ot Q(j)$, $\we^3Q(j)$, $Q\ot\we^2Q(j)$, $Q\ot\we^3Q(j)$, $\we^2Q\ot\we^2Q(j)$,
$\we^2Q\ot\we^3Q(j)$ for $j=0,\ldots,4$;

\noindent
(iii) $Q^{\ast}(j)$, $Q^{\ast}\ot Q(j)$, $Q^{\ast}\ot\we^2 Q(j)$ for $j=1,\ldots, 5$;

\noindent
(iv) $Q\ot Q\ot Q$, $Q\ot Q\ot\we^2Q$, $Q\ot Q\ot\we^3Q$, $Q^{\ast}\ot Q\ot Q(1)$, 
$Q\ot\we^2Q\ot\we^2Q$, $Q\ot\we^2Q\ot\we^3Q$, $Q^{\ast}\ot Q\ot\we^2Q(1)$.

The fact that all these objects belong to $\CC$ follows from Lemmas \ref{L5-gen-1}, 
\ref{L5-gen-4}--\ref{L5-gen-8} below.
\ed

In the following Lemmas we often use the fact that $\CC$ is closed under direct
summands (as an admissible subcategory). Also, by a resolution of $S^nQ$ we mean
the exact sequence
$$\ldots\to\we^2Q^{\ast}\ot S^{n-2}V\ot\OO\to Q^{\ast}\ot S^{n-1}V\ot\OO\to S^nV\ot\OO\to S^nQ\to 0.$$
By the standard filtration of $\we^k(V\ot\OO)$ we mean the filtration associated with exact sequence
\eqref{basic-seq}. This filtration has vector bundles $\we^iQ^{\ast}\ot\we^{k-i}Q$ as consecutive quotients.
Recall also that $\we^5Q=\OO(1)$, so we have isomorphisms 
$\we^iQ^{\ast}(1)\simeq\we^{5-i}Q$.

\begin{lem}\label{L5-gen-1} 
(i) For $j=0,\ldots, 5$ the following objects are in $\CC$: $\OO(j)$, $Q(j)$, $\we^2Q(j)$,
$\we^3Q(j)$, $Q\ot\we^2Q(j)$, 
$Q^{\ast}(j)$, $Q^{\ast}\ot\we^2Q(j)$, $S^2Q^{\ast}(j)$.

\noindent
(ii) For $j=1,\ldots,5$ the following objects are in $\CC$: $Q^{\ast}\ot Q^{\ast}(j)$, 
$Q^{\ast}\ot Q(j)$, $Q\ot Q(j)$, $S^nQ(j)$ for $n\ge 2$.

\noindent
(iii) For $j=0,\ldots,4$ one has $Q\ot\we^3Q(j)\in\CC$ and $Q^{\ast}\ot\we^3Q(j)\in\CC$.

\noindent
(iv) For $j=1,\ldots,4$ one has 
$\we^3Q\ot\we^2Q(j-1)=\we^2Q^{\ast}\ot\we^2Q(j)\in\CC$ and $S^2Q\ot\we^2Q(j)\in\CC$.

\noindent
(v) For $j=1,\ldots,5$ and for $n\ge 2$ one has $Q\ot S^nQ(j)\in\CC$ and $Q^{\ast}\ot S^nQ(j)\in\CC$.

\noindent
(vi) For $j=1,\ldots,5$ the following objects are in $\CC$: 
$Q\ot Q\ot Q(j)$, $Q^{\ast}\ot Q\ot Q(j)$, $Q^{\ast}\ot Q^{\ast}\ot Q(j)$ and $Q^{\ast}\ot Q^{\ast}\ot Q^{\ast}(j)$.
\end{lem}

\Pf . (i) To check the assertion for $Q\ot\we^2Q(j)$ we 
observe that $R_1=S^{2,1}Q$ is contained in 
$\lan Q,F_1\ran$ as follows from exact
sequence \eqref{FR-seq}. This implies that
$Q\ot\we^2Q(j)=\we^3Q(j)\oplus S^{2,1}Q(j)$ belongs to $\CC$ for $j=0,\ldots,5$.

The assertions for $Q^{\ast}(j)$ and $Q^{\ast}\ot\we^2Q(j)$ follow from the sequence \eqref{basic-seq}.
The assertion for $S^2Q^{\ast}(j)$ follows by considering the dual sequence to the resolution of $S^2Q$.

\noindent
(ii) Use the decomposition $Q^{\ast}\ot Q^{\ast}(j)=S^2Q^{\ast}(j)\oplus\we^2Q^{\ast}(j)=
S^2Q^{\ast}(j)\oplus\we^3Q(j-1)$ and (i). Then use sequence \eqref{basic-seq}.
For $S^nQ(j)$ the assertion is checked using part (i) and the resolution of $S^nQ$.

\noindent
(iii) To prove the assertion for $Q\ot\we^3Q(j)$ 
use the isomorphism $Q\ot\we^3Q(j)\equiv Q\ot\we^2Q^{\ast}(j+1)$ and consider the standard filtration of
$\we^3(V\ot\OO)$ tensored with $\OO(j+1)$ (and then use part (i)).
For the second assertion use sequence \ref{basic-seq}.

\noindent
(iv) To check that $\we^2Q^{\ast}\ot\we^2Q(j)\in\CC$ 
use the standard filtration of $\we^4(V\ot\OO)$ tensored with $\OO(j)$. Next, to derive that 
$S^2Q\ot\we^2Q(j)\in\CC$
use resolution of $S^2Q$.

\noindent
(v) For $Q\ot S^nQ(j)$ use the resolution for $S^nQ$   tensored with $Q(j)$ and parts
(i), (ii) and (iii). For $Q^{\ast}\ot S^nQ(j)$ use sequence \eqref{basic-seq} and part (ii).

\noindent
(vi) The assertion for $Q\ot Q\ot Q(j)$ follows from
 the decomposition $Q\ot Q\ot Q(j)=Q\ot \we^2Q(j)\oplus Q\ot S^2Q(j)$
and parts (i) and (v). The rest follows using sequence \eqref{basic-seq} and part (ii).
\ed

\begin{lem}\label{L5-gen-2}
(i) One has $S^{3,1,1}Q\in\CC$ and $S^{3,1,1}(3)Q\in\CC$.

\noindent
(ii) One has $S^2Q\ot\we^3Q\in\CC$ and $S^2Q\ot\we^3Q(3)\in\CC$.

\noindent
(iii) One has $\we^2 Q^{\ast}\ot\we^3Q\in\CC$ and $\we^2 Q^{\ast}\ot\we^3Q(3)\in\CC$.
\end{lem}

\Pf . (i) First, exact sequence \eqref{PG-seq} shows that $P,P(3)\in\CC$. Next,
exact sequence \eqref{PT-seq} shows that $T,T(3)\in\CC$, where
$T=S^{3,1,1}Q$.

\noindent
(ii) Since we have the decomposition 
$$S^2Q\ot\we^3Q=S^{3,1,1}Q\oplus S^{2,1,1,1}Q,$$ 
part (i) shows that it is enough to check the similar assertion for $S^{2,1,1,1}Q$.
But $S^{2,1,1,1}Q$ is a direct summand in $Q\ot\we^4Q=Q\ot Q^{\ast}(1)$, so
the statement follows from Lemma \ref{L5-gen-1}(ii).

\noindent
(iii) This follows from (ii) using resolution for $S^2Q$ and Lemma \ref{L5-gen-1}(iii).
\ed

\begin{lem}\label{L5-gen-3}
(i) For $j=1,2,3$, $\we^3Q\ot\we^3Q(j)\in\CC$ if and only if $\we^2Q\ot\we^2Q(j)\in\CC$.

\noindent
(ii) For $j=1,\ldots,4$, $\we^2Q\ot\we^2Q(j)\in\CC$ if and only if $S^2Q\ot S^2Q(j)\in\CC$.

\noindent
(iii) For $j=1,\ldots,4$, the following conditions are equivalent:
\begin{enumerate}
\item[(1)]  $\we^2Q^{\ast}\ot\we^3Q(j-1)\in\CC$; 
\item[(2)] $\we^2Q^{\ast}\ot S^2Q(j)\in\CC$;
\item[(3)] $Q\ot Q\ot\we^2Q(j)\in\CC$. 
\end{enumerate}
\end{lem}

\Pf .
(i) Use the standard filtration of $\we^5(V\ot\OO)$ tensored with $\OO(j+1)$ and Lemma \ref{L5-gen-1}(ii).

\noindent
(ii) Use the decompositions 
$$S^2Q\ot S^2Q=Q\ot S^3Q\oplus S^{2,2}Q,\ \ \we^2Q\ot\we^2Q=Q\ot\we^3Q\oplus S^{2,2}Q$$
and Lemma \ref{L5-gen-1}(iii),(v).

\noindent
(iii) First, the equivalence of (1) and (2) follows by considering
the resolution of $S^2Q$ and using Lemma \ref{L5-gen-1}(iii). Next,
we observe that 
$$\we^2Q^{\ast}\ot Q\ot Q(j)=\we^2Q^{\ast}\ot\we^2Q(j)\oplus\we^2Q^{\ast}\ot S^2Q(j)$$
and that $\we^2Q^{\ast}\ot\we^2Q(j)\in\CC$ for $j=1,\ldots,4$ by Lemma \ref{L5-gen-1}(iv).
Therefore, (2) is equivalent to $\we^2Q^{\ast}\ot Q\ot Q(j)\in\CC$.
On the other hand, sequence \eqref{basic-seq} and Lemma \ref{L5-gen-1}(i) imply that in condition (3)
we can replace $Q\ot Q\ot\we^2Q(j)$ with $Q^{\ast}\ot Q\ot\we^2Q(j)$.
Now the equivalence of (2) and (3) follows by considering
the standard filtration of $\we^3(V\ot\OO)$ tensored with $Q(j)$ and using Lemma \ref{L5-gen-1}(i),(iii).
\ed

\begin{lem}\label{L5-gen-4}
(i) For $j=0,1,2,3$ one has $\we^3Q\ot\we^3Q(j-1)\in\CC$, $S^2Q\ot\we^3Q(j)\in\CC$ and
$S^2Q\ot S^2Q(j+1)\in\CC$.

\noindent
(ii) For $j=1,2,3,4$ one has $\we^2Q\ot Q\ot Q(j)\in\CC$ and $\we^2Q\ot Q^{\ast}\ot Q(j)\in\CC$.

\noindent
(iii) One has $\we^2Q\ot\we^2Q(j)\in\CC$, $\we^3Q\ot\we^3Q(j-1)\in\CC$ for $j=1,\ldots,4$, and
$S^2Q\ot\we^3Q(j)\in\CC$ for $j=0,\ldots,4$.

\noindent
(iv) One has $\we^3Q\ot Q\ot Q(j)\in\CC$ and $\we^3Q\ot Q^{\ast}\ot Q(j)\in\CC$
for $j=0,1,2,3$. 

\noindent
(v) One has $\we^3Q\ot\we^2Q\ot Q(j)\in\CC$ for $j=1,2$.
\end{lem}

\Pf . (i) For $j=0$ and $j=3$ the first assertion follows from Lemma \ref{L5-gen-2}(iii). 
By Lemma \ref{L5-gen-3}(iii) this implies
that $S^2Q\ot\we^3Q\in\CC$ and $S^2Q\ot\we^3Q(3)\in\CC$. Next, using
the resolution for $S^2Q$ and Lemma \ref{L5-gen-1}(v)
we obtain $S^2Q\ot S^2Q(1)\in\CC$ and $S^2Q\ot S^2Q(4)\in\CC$.
By Lemma \ref{L5-gen-3}(ii), this implies that $\we^2Q\ot\we^2Q(1)\in\CC$ and
$\we^2Q\ot\we^2Q(4)\in\CC$. By Lemma \ref{L5-gen-3}(i), it follows that $\we^3Q\ot\we^3Q(1)\in\CC$,
which also leads to $S^2Q\ot\we^3Q(2)\in\CC$ and $S^2Q\ot S^2Q(3)\in\CC$ as before.

On the other hand, combining Lemma \ref{L5-gen-3}(i) with Lemma \ref{L5-gen-2}(iii)
we also get $\we^2Q\ot\we^2Q(2)\in\CC$.
By Lemma \ref{L5-gen-3}(ii), this implies that $S^2Q\ot S^2Q(2)\in\CC$. Considering the resolution
for $S^2Q$ this leads to $S^2Q\ot\we^3Q(1)\in\CC$ and $\we^3Q\ot\we^3Q\in\CC$ as before.

\noindent
(ii) The first assertion immediately follows from (i) and from Lemma \ref{L5-gen-3}(iii). The second
follows from the first using sequence \eqref{basic-seq}.

\noindent
(iii) This follows from (i), (ii) and Lemma \ref{L5-gen-3}(i).

\noindent
(iv) Using sequence \eqref{basic-seq} and Lemma \ref{L5-gen-1}(iii) we see that it is enough
to show that $\we^3Q\ot Q\ot Q(j)\in\CC$. To this end we use the decomposition
$$\we^3Q\ot Q\ot Q(j)=\we^3Q\ot\we^2Q(j)\oplus\we^3Q\ot S^2Q(j),$$
part (iii) and Lemma \ref{L5-gen-1}(iv).

\noindent
(v) We start with the isomorphism $\we^3Q\ot\we^2Q\ot Q(j)\simeq\we^2Q^{\ast}\ot\we^2Q\ot Q(j+1)$.
Now the assertion follows by considering the standard filtration of $\we^4(V\ot\OO)$ 
tensored with $Q(j+1)$ and using parts (ii), (iv) and
Lemma \ref{L5-gen-3}(ii). 
\ed

\begin{lem}\label{L5-gen-5}
(i) For $j=1,2,3$ one has $\we^2Q\ot\we^2Q\ot Q(j)\in\CC$.

\noindent
(ii) One has $\we^2Q\ot\we^2Q\ot \we^2Q(2)\in\CC$.

\noindent
(iii) One has $\we^2Q\ot\we^2Q\ot S^2Q(2)\in\CC$.

\noindent
(iv) One has $\we^3Q\ot\we^2Q(4)\in\CC$.
\end{lem}

\Pf .(i) Suppose first that $j=1$.
Then considering the filtration of $\we^5(V\ot\OO)\ot Q(2)$ and using Lemma \ref{L5-gen-1}(vi), as well
as the fact that $Q^{\ast}\ot Q^{\ast}\ot\we^2Q(3)\in\CC$ (which is a consequence of Lemma
\ref{L5-gen-4}(ii)), we reduce ourselves to showing that $\we^3Q\ot\we^3Q\ot Q(1)\in\CC$.
Now the isomorphism $\we^3Q\ot\we^3Q\ot Q(1)\simeq \we^3Q\ot\we^2Q^{\ast}\ot Q(2)$
and the standard filtration of $\we^3Q\ot\we^3(V\ot\OO)(2)$ show that it is enough to check that the following
objects are in $\CC$:
$$\we^3Q\ot\we^2Q\ot Q^{\ast}(2),\ \we^3Q\ot\we^3Q(2),\ \we^3Q\ot\we^3Q^{\ast}(2).$$
For the second and the third this follows from Lemma \ref{L5-gen-4}(iii) and Lemma \ref{L5-gen-1}(iv),
respectively.
For the first object this follows from Lemmas \ref{L5-gen-1}(iv) and \ref{L5-gen-4}(v) using \eqref{basic-seq}.

Now consider the case $j=2$ or $j=3$. 
By sequence \eqref{basic-seq} and Lemma \ref{L5-gen-4}(iii), it is enough to
prove that $\we^2Q\ot\we^2Q\ot Q^{\ast}(j)\in\CC$. Now the standard filtration of $\we^2Q\ot\we^3(V\ot\OO)(j)$
shows that it is enough to check that the following objects are in $\CC$:
$$\we^2Q\ot\we^3Q(j),\ \we^2Q\ot\we^3Q^{\ast}(j),\ \we^2Q\ot\we^2Q^{\ast}\ot Q(j).$$
But this follows from Lemmas \ref{L5-gen-1}(iv), \ref{L5-gen-4}(iii) and \ref{L5-gen-4}(v), respectively.

\noindent
(ii) First, considering the standard filtration of $\we^5(V\ot\OO)\ot\we^2Q(3)$, 
we reduce ourselves to showing that the following objects are in $\CC$:
$$\we^3Q\ot\we^3Q\ot\we^2Q(2),\ Q\ot Q\ot\we^2Q(2),\ Q^{\ast}\ot Q^{\ast}\ot\we^2Q(4).$$
For the second and the third this follows from Lemma \ref{L5-gen-4}(ii). Now using the isomorphism
$\we^3Q\ot\we^3Q\ot\we^2Q(2)\simeq\we^3Q\ot\we^2Q^{\ast}\ot\we^2Q(3)$ and the standard filtration of
$\we^3Q\ot\we^4(V\ot\OO)(3)$ we are led to showing that the following objects are in $\CC$:
$$\we^3Q\ot\we^3Q\ot Q^{\ast}(3),\ \we^3Q\ot\we^2Q\ot Q(2),\ \we^3Q\ot Q(2),\ \we^3Q\ot Q^{\ast}(4).$$
For the second object this follows from Lemma \ref{L5-gen-4}(v), while for the last two it follows
from Lemma \ref{L5-gen-1}(iii). Thus, it remains to check that $\we^3Q\ot\we^3Q\ot Q^{\ast}(3)\in\CC$.
Using the standard filtration of $\we^5(V\ot\OO)\ot Q^{\ast}(4)$ we see
that it is enough to verify that the following objects are in $\CC$:
$$\we^2Q\ot\we^2Q\ot Q^{\ast}(3),\ Q\ot Q\ot Q^{\ast}(3),\ Q^{\ast}(3), Q^{\ast}(5),\
Q^{\ast}\ot Q^{\ast}\ot Q^{\ast}(5).$$
For the second and the last object this follows from Lemma \ref{L5-gen-1}(vi).
On the other hand, using \eqref{basic-seq}, part (i) and Lemma \ref{L5-gen-4}(iii) 
we see that $\we^2Q\ot\we^2Q\ot Q^{\ast}(3)\in\CC$.

\noindent
(iii) First, using the resolution for $S^2Q$ we reduce the problem to
showing that $\we^2Q\ot\we^2Q\ot\we^2Q^{\ast}(2)\in\CC$ (here we also use 
part (i), sequence \eqref{basic-seq} and Lemma \ref{L5-gen-4}(iii)). 
Next, the standard filtration of $\we^2Q\ot\we^4(V\ot\OO)(2)$ shows that it is enough to check that
the following objects are in $\CC$:
$$\we^2Q\ot\we^3Q\ot Q^{\ast}(2),\ \we^2Q\ot \we^2Q\ot Q(1),\ \we^2Q\ot Q^{\ast}(3),\ \we^2Q\ot Q(1).$$
For the last two objects this follows from Lemma \ref{L5-gen-1}(i).
For the second object the assertion follows from part (i).
Finally, to check that $\we^2Q\ot\we^3Q\ot Q^{\ast}(2)\in\CC$ we use 
sequence \eqref{basic-seq}, Lemma \ref{L5-gen-4}(v) and Lemma \ref{L5-gen-1}(iv).

\noindent
(iv) Let us start with the decomposition 
$$\we^3Q\ot\we^2Q(4)=\OO(5)\oplus S^{2,1,1,1}Q(4)\oplus S^{2,2,1}Q(4).$$
Now observe that $S^{2,1,1,1}Q(4)$ is a direct summand in $Q\ot\we^4Q(4)=Q\ot Q^{\ast}(5)$ which
is in $\CC$ by Lemma \ref{L5-gen-1}(ii), while
$S^{2,2,1}Q(4)$ is a direct summand in $S^2Q\ot S^2Q\ot Q(4)$. Using the resolution of $S^2Q$
we reduce ourselves to checking that the following objects are in $\CC$:
$$S^2Q\ot\we^2Q^{\ast}\ot Q(4),\ S^2Q\ot Q(4),\ S^2Q\ot Q^{\ast}\ot Q(4).$$
For the second object this follows from Lemma \ref{L5-gen-1}(v).
Using \eqref{basic-seq} we can replace the third object by 
$S^2Q\ot Q\ot Q(4)=S^2Q\ot\we^2Q(4)\oplus S^2Q\ot S^2Q(4)$ which is in $\CC$
by Lemmas \ref{L5-gen-1}(iv) and \ref{L5-gen-4}(i).
Next, we use the isomorphism
$S^2Q\ot\we^2Q^{\ast}\ot Q(4)\simeq S^2Q\ot\we^3Q\ot Q(3)$ and the resolution of $S^2Q$
to reduce the problem to showing that the following objects are in $\CC$:
$$\we^2Q^{\ast}\ot\we^3Q\ot Q(3),\ \we^3Q\ot Q(3),\ \we^3Q\ot Q^{\ast}\ot Q(3).$$
The second and third objects are in $\CC$ by Lemmas \ref{L5-gen-1}(iii) and
\ref{L5-gen-4}(iv), respectively.
For the first object we use the isomorphism
$\we^2Q^{\ast}\ot\we^3Q\ot Q(3)\simeq \we^3Q\ot\we^3Q\ot Q(2)$ and the standard filtration of
$\we^5(V\ot\OO)\ot Q(3)$ to reduce ourselves to proving that the following objects are in $\CC$:
$$\we^2Q\ot\we^2Q\ot Q(2),\ Q\ot Q\ot Q(2),\ Q^{\ast}\ot Q^{\ast}\ot Q(4).$$
For the first object this follows from (i), and for the second and the third---from
Lemma \ref{L5-gen-1}(vi).
\ed

\begin{lem}\label{L5-gen-6}
(i) One has $S^3Q\ot S^3Q(2)\in\CC$.

\noindent
(ii) One has $\we^2Q\ot\we^2Q\in\CC$.

\noindent
(iii) One has $Q\ot Q\in\CC$.
\end{lem}

\Pf . (i) Consider the decomposition 
$$S^3Q\ot S^3Q(2)=S^6Q(2)\oplus S^{5,1}Q(2)\oplus S^{4,2}Q(2)\oplus S^{3,3}Q(2).$$
By Lemma \ref{L5-gen-1}(ii), we have $S^6Q(2)\in\CC$.
Next, we observe that $S^{5,1}Q(2)$ is a direct summand in $S^4Q\ot\we^2Q(2)$ and
use the resolution of $S^4Q$ to deduce that this object is in $\CC$ from the inclusions
$Q\ot\we^2Q(1)\in\CC$, $Q^{\ast}\ot\we^2Q(2)\in\CC$, $\we^2Q^{\ast}\ot\we^2Q(2)\in\CC$, $\we^3Q^{\ast}\ot\we^2Q(2)\in\CC$, that follow from Lemmas \ref{L5-gen-1}(i), \ref{L5-gen-1}(iv) and
\ref{L5-gen-4}(iii). Finally, we note that $S^{4,2}Q(2)\oplus S^{3,3}Q(2)$ is a direct summand in
$$\we^2Q\ot\we^2Q\ot Q\ot Q(2)=\we^2Q\ot\we^2Q\ot\we^2Q(2)\oplus\we^2Q\ot\we^2Q\ot S^2Q(2)$$
which is in $\CC$ by Lemma \ref{L5-gen-5}(ii),(iii).

\noindent
(ii) We use the isomorphism $\we^2Q\ot\we^2Q\simeq\we^3Q^{\ast}\ot\we^3Q^{\ast}(2)$
and then use the resolution of $S^3Q$ twice to relate this to $S^3Q\ot S^3Q(2)$ which is in $\CC$
by part (i). It remains to check that the objects
that appear in between, namely,
$$S^3Q(2),\ S^3Q\ot Q^{\ast}(2),\ S^3Q\ot\we^2Q^{\ast}(2),\  
\we^3Q^{\ast}(2),\ Q^{\ast}\ot\we^3Q^{\ast}(2),
\ \we^2Q^{\ast}\ot\we^3Q^{\ast}(2),$$
are all in $\CC$. For the last three objects this follows from Lemma \ref{L5-gen-1}(i),(iv), while
for the first three one has to use the resolution of $S^3Q$ to reduce to the objects we have already
dealt with.

\noindent
(iii) The standard filtration of $\we^5(V\ot\OO)(1)$ reduces the problem to showing that 
$\we^2Q\ot\we^2Q$ and $\we^3Q\ot\we^3Q$
are in $\CC$ (where we also use Lemma \ref{L5-gen-1}(ii)). It remains to apply 
part (ii) and Lemma \ref{L5-gen-4}(iii).
\ed

\begin{lem}\label{L5-gen-7}
(i) One has $S^3Q\ot S^2Q(1)\in\CC$.

\noindent
(ii) One has $\we^2Q\ot S^2Q\in\CC$.

\noindent
(iii) One has $\we^2Q\ot Q\ot Q\in\CC$.

\noindent
(iv) One has $\we^3Q\ot\we^2Q\ot Q\in\CC$.
\end{lem}

\Pf . (i) Consider the decomposition
$$S^3Q\ot S^2Q(1)=S^5Q(1)\oplus S^{4,1}Q(1)\oplus S^{3,2}Q(1).$$
By Lemma \ref{L5-gen-1}(ii), we have $S^5Q(1)\in\CC$.
On the other hand, $S^{4,1}Q(1)$ is a direct summand in $S^3Q\ot\we^2Q(1)$.
The resolution of $S^3Q$ relates the latter object to $Q^{\ast}\ot\we^2Q(1)$,
$\we^2Q^{\ast}\ot\we^2Q(1)$ and $\we^3Q^{\ast}\ot\we^2Q(1)$ which are all
 in $\CC$ (for the last one use Lemma \ref{L5-gen-6}(ii)). Finally,
 $S^{3,2}Q(1)$ is a direct summand in $\we^2Q\ot\we^2Q\ot Q(1)$ which is in $\CC$ by
 Lemma \ref{L5-gen-5}(i).

\noindent
(ii) Using the resolution for $S^3Q$ we can relate $\we^2Q\ot S^2Q=\we^3Q^{\ast}\ot S^2Q(1)$
with $S^3Q\ot S^2Q(1)$, which is in $\CC$ by part (i). The objects appearing in between,
namely, $Q^{\ast}\ot S^2Q(1)$ and $\we^2Q^{\ast}\ot S^2Q(1)$ are in $\CC$,
by Lemmas \ref{L5-gen-1}(vi), \ref{L5-gen-2}(ii).

\noindent
(iii) Since $\we^2Q\ot Q\ot Q=\we^2Q\ot\we^2Q\oplus\we^2Q\ot S^2Q$, this follows from
part (ii) and Lemma \ref{L5-gen-6}(ii).

\noindent
(iv) Considering the filtration of $\we^4(V\ot\OO)\ot Q(1)$ we reduce ourselves to showing that
the following objects are in $\CC$:
$$\we^2Q\ot Q\ot Q,\ Q\ot Q,\ Q^{\ast}\ot Q(2),\ Q^{\ast}\ot Q\ot\we^3Q(1).$$
Now the assertion follows from part (iii) and Lemmas \ref{L5-gen-6}(iii), \ref{L5-gen-1}(ii) and
\ref{L5-gen-4}(iv).
\ed

\begin{lem}\label{L5-gen-8}
(i) One has $S^2Q\ot S^4Q(1)\in\CC$.

\noindent
(ii) One has $S^2Q\ot Q\in\CC$.

\noindent
(iii) One has $Q\ot Q\ot Q\in\CC$.

\noindent
(iv) One has $\we^2Q\ot\we^2Q\ot Q\in\CC$.
\end{lem}

\Pf . (i) Consider the decomposition 
$$S^2Q\ot S^4Q(1)=S^6Q(1)\oplus S^{5,1}Q(1).$$
By Lemma \ref{L5-gen-1}(ii), we have $S^6Q(1)\in\CC$. On the other hand,
$S^{5,1}Q(1)$ is a direct summand in $\we^2Q\ot S^4Q(1)$. Using
the resolution of $S^4Q$ we reduce the problem to checking that
the following objects are in $\CC$:
$$Q\ot\we^2Q,\ \we^2Q\ot\we^2Q,\ \we^2Q^{\ast}\ot\we^2Q(1),\ Q^{\ast}\ot\we^2Q(1),\
\we^2Q(1),$$
which follows from our previous work (for the second object use Lemma \ref{L5-gen-6}(ii)).

\noindent
(ii) Tensoring the resolution for $S^4Q$ with $S^2Q(1)$ we get an exact sequence
\begin{align*}
&0\to S^2Q\ot Q\to V\ot S^2Q\ot\we^2Q\to S^2V\ot S^2Q\ot\we^3Q\to S^3V\ot S^2Q\ot Q^{\ast}(1)\to \\
&S^4V\ot S^2Q(1)\to S^2Q\ot S^4Q(1)\to 0
\end{align*}
By part (i), one has $S^2Q\ot S^4Q(1)\in\CC$. Next, $S^2Q\ot Q^{\ast}(1)$ and $S^2Q(1)$ are in $\CC$
by Lemma \ref{L5-gen-1}(v),(ii). Finally, $S^2Q\ot\we^2Q$ and $S^2Q\ot\we^3Q$ 
are in $\CC$ by Lemmas \ref{L5-gen-7}(ii) and \ref{L5-gen-4}(iii). 
Hence, $S^2Q\ot Q\in\CC$.

\noindent
(iii) This follows from the decomposition 
$$Q\ot Q\ot Q=S^2Q\ot Q\oplus \we^2Q\ot Q,$$
part (ii) and Lemma \ref{L5-gen-1}(i).

\noindent
(iv) This is proved by the same method as the case $j=1$ of Lemma \ref{L5-gen-5}(i),
using part (iii).
\ed

\end{document}